**Some possible approaches to the Riemann Hypothesis via the Li/Keiper constants**

Donal F. Connon

dconnon@btopenworld.com

18 February 2010

**Abstract**

In this paper we consider some possible approaches to the proof of the Riemann Hypothesis using the Li criterion.

Some examples of the potentially useful formulae are set out below:

$$\lambda_n = \frac{n}{2}(\gamma + \log \pi) + \sum_{m=2}^{n} \binom{n}{m} 2^{-m} \varsigma(m) + \sum_{m=1}^{n} (-1)^m \binom{n}{m} \frac{1}{(m-1)!} \frac{d^m}{ds^m} \log[(s-1)\varsigma(s)]\Big|_{s=0}$$

$$\lambda_n = n\lambda_1 + \sum_{m=2}^{n} \binom{n}{m}[2^{-m}\varsigma(m) + (-1)^m b_{m-1}]$$

where $\lambda_n$ and $b_n$ are the Li/Keiper constants and the Lehmer constants respectively.

CONTENTS                                                                 Page

**1. Introduction**

The Riemann xi function $\xi(s)$ is defined as

$$(1.1) \qquad \xi(s) = \frac{1}{2} s(s-1)\pi^{-s/2} \Gamma(s/2)\varsigma(s)$$

and we see from the functional equation for the Riemann zeta function $\varsigma(s)$

$$(1.2) \qquad \varsigma(1-s) = 2(2\pi)^{-s}\Gamma(s)\cos(\pi s/2)\varsigma(s)$$

that $\xi(s)$ satisfies the functional equation

$$(1.3) \qquad \xi(s) = \xi(1-s)$$

In 1996, Li [29] defined the sequence of numbers $(\lambda_n)$ by

$$(1.4) \qquad \lambda_n = \frac{1}{(n-1)!}\frac{d^n}{ds^n}[s^{n-1}\log\xi(s)]\Big|_{s=1}$$

and proved that a necessary and sufficient condition for the non-trivial zeros $\rho$ of the Riemann zeta function to lie on the critical line $s = \frac{1}{2} + i\tau$ is that $\lambda_n$ is non-negative for every positive integer $n$. Earlier in 1991, Keiper [25] showed that if the Riemann hypothesis is true, then $\lambda_n > 0$ for all $n \geq 1$.

Using the Leibniz rule for differentiation we may write (1.4) as

$$(1.4.1) \qquad \lambda_n = \sum_{m=1}^{n}\binom{n}{m}\frac{1}{(m-1)!}\frac{d^m}{ds^m}\log\xi(s)\Big|_{s=1}$$

Li also showed that

$$(1.5) \qquad \lambda_n = \sum_{\rho}\left[1 - \left(1 - \frac{1}{\rho}\right)^n\right]$$

Taking logarithms of (1.1) and noting that $\frac{s}{2}\Gamma\left(\frac{s}{2}\right) = \Gamma\left(1 + \frac{s}{2}\right)$ we see that

$$(1.6) \qquad \log\xi(s) = \log\Gamma\left(1 + \frac{s}{2}\right) - \frac{s}{2}\log\pi + \log[(s-1)\varsigma(s)]$$

We have the Maclaurin expansion about $s = 1$

$$(1.7.1) \qquad \log\xi(s) = -\log 2 - \sum_{k=1}^{\infty}(-1)^k\frac{\sigma_k}{k}(s-1)^k$$

and about $s = 0$

(1.7.2)     $\log \xi(s) = -\log 2 - \sum_{k=1}^{\infty} \frac{\sigma_k}{k} s^k$

where the constant term in (1.7) arises because

$$\lim_{s \to 1} \xi(s) = \frac{1}{2} \pi^{-1/2} \Gamma(1/2) \lim_{s \to 1}[(s-1)\varsigma(s)] = \frac{1}{2}$$

We note that the coefficients $\sigma_k$ are defined by

$$\sigma_k = \frac{(-1)^{k+1}}{(k-1)!} \frac{d^k}{ds^k} \log \xi(s) \bigg|_{s=1}$$

and we see that

$$\frac{d}{ds} \log \xi(s) = \frac{\xi'(s)}{\xi(s)} = -\sum_{k=1}^{\infty}(-1)^k \sigma_k (s-1)^{k-1}$$

and comparing this with (1.4) we immediately see that $\lambda_1 = \sigma_1$.

We see that $\frac{\xi'(0)}{\xi(0)} = \sum_{k=1}^{\infty} \sigma_k$ and, since the series is convergent, we deduce that

$$\lim_{k \to \infty} \sigma_k = 0$$

Having regard to the definition of the Li/Keiper constants

$$\lambda_n = \frac{1}{(n-1)!} \frac{d^n}{ds^n}[s^{n-1} \log \xi(s)] \bigg|_{s=1}$$

we consider

$$s^{n-1} \log \xi(s) = -s^{n-1} \log 2 - \sum_{k=1}^{\infty} \frac{\sigma_k}{k} s^{k+n-1} = -s^{n-1} \log 2 - \sum_{k=1}^{\infty}(-1)^k \frac{\sigma_k}{k} s^{n-1}(s-1)^k$$

and using the Leibniz differentiation formula we immediately obtain

(1.8)                $\lambda_n = -\sum_{m=1}^{n}(-1)^m \binom{n}{m} \sigma_m$

We have from Keiper's paper [25]

(1.9)
$$\sigma_n = \sum_\rho \frac{1}{\rho^n}$$

and since

$$\lambda_n = \sum_\rho \left[ 1 - \left( 1 - \frac{1}{\rho} \right)^n \right] = -\sum_\rho \sum_{k=0}^n \binom{n}{k} \frac{(-1)^k}{\rho^k}$$

$$= -\sum_{k=0}^n \binom{n}{k} (-1)^k \sum_\rho \frac{1}{\rho^k}$$

we deduce (1.8) again and we note that this also implies that $\lambda_1 = \sigma_1$.

The eta constants $\eta_k$ are defined by reference to the logarithmic derivative of the Riemann zeta function

(1.10)
$$\frac{d}{ds}[\log \varsigma(s)] = \frac{\varsigma'(s)}{\varsigma(s)} = -\frac{1}{s-1} - \sum_{k=1}^\infty \eta_{k-1}(s-1)^{k-1} \qquad |s-1| < 3$$

and we may also note that this is equivalent to

(1.11)
$$\frac{d}{ds}\log[(s-1)\varsigma(s)] = \frac{\varsigma'(s)}{\varsigma(s)} + \frac{1}{s-1} = -\sum_{k=1}^\infty \eta_{k-1}(s-1)^{k-1}$$

and, noting that $\lim_{s \to 1}[(s-1)\varsigma(s)] = 1$, we obtain upon integration

(1.12)
$$\log[(s-1)\varsigma(s)] = -\sum_{k=1}^\infty \frac{\eta_{k-1}}{k}(s-1)^k$$

We also have the Hadamard representation of the Riemann zeta function

(1.13)
$$\varsigma(s) = \frac{e^{cs}}{2(s-1)\Gamma\left(1+\frac{s}{2}\right)} \prod_\rho \left(1 - \frac{s}{\rho}\right) e^{s/\rho}$$

where $c = \log(2\pi) - 1 - \gamma/2$.

It was shown by Zhang and Williams [38] in 1994 that (as corrected by Coffey [13])

(1.14) $\frac{d}{ds}\log[(s-1)\varsigma(s)] = \frac{\varsigma'(s)}{\varsigma(s)} + \frac{1}{s-1} = \gamma + \sum_{n=1}^\infty (-1)^n \left[ \left(1 - \frac{1}{2^{n+1}}\right)\varsigma(n+1) + \sigma_{n+1} - 1 \right](s-1)^n$

and we see that

$$\lim_{s \to 1} \frac{d^{n+1}}{ds^{n+1}} \log[(s-1)\varsigma(s)] = (-1)^n n! \left[ \sigma_{n+1} + \left(1 - \frac{1}{2^{n+1}}\right)\varsigma(n+1) - 1 \right]$$

Differentiating (1.11) gives us

$$\lim_{s \to 1} \frac{d^{n+1}}{ds^{n+1}} \log[(s-1)\varsigma(s)] = -n!\eta_n$$

and hence we have for $n \geq 1$

(1.15)       $$\eta_n = (-1)^{n+1}\left[ \sigma_{n+1} + \left(1 - \frac{1}{2^{n+1}}\right)\varsigma(n+1) - 1 \right]$$

We have the representation for $\lambda_n$ given by Bombieri and Lagarias [6] in 1999 and by Coffey ([11] and [12]) in 2004

(1.16)       $$\lambda_n = 1 - \frac{n}{2}[\log \pi + \gamma + 2\log 2] + S_1(n) + S_2(n)$$

where

(1.16.1)     $$S_1(n) = \sum_{m=2}^{n}\binom{n}{m}(-1)^m\left(1 - \frac{1}{2^m}\right)\varsigma(m)$$

(1.16.2)     $$S_2(n) = \sum_{m=1}^{n}\binom{n}{m}\frac{1}{(m-1)!}\frac{d^m}{ds^m}\log[(s-1)\varsigma(s)]\Big|_{s=1}$$

Maślanka [30] makes the decomposition

$$\lambda_m = \bar{\lambda}_m + \widetilde{\lambda}_m$$

where the "trend" $\bar{\lambda}_m$ is strictly increasing

$$\bar{\lambda}_m = 1 - \frac{m}{2}[\log \pi + \gamma + 2\log 2] + S_1(m)$$

and the "oscillations" $\widetilde{\lambda}_m$ are given by

$$\widetilde{\lambda}_m = -\sum_{n=1}^{m}\binom{m}{n}\eta_{n-1}$$

From (1.11) we see that

$$\lim_{s\to 1}\frac{1}{(n-1)!}\frac{d^n}{ds^n}\log[(s-1)\varsigma(s)] = -\eta_{n-1}$$

and therefore we have

$$-\sum_{n=1}^{m}\binom{m}{n}\eta_{n-1} = \lim_{s\to 1}\sum_{n=1}^{m}\binom{m}{n}\frac{1}{(n-1)!}\frac{d^n}{ds^n}\log[(s-1)\varsigma(s)]$$

This then shows that

(1.17) $$\widetilde{\lambda}_m = S_2(m) = -\sum_{n=1}^{m}\binom{m}{n}\eta_{n-1}$$

The generalised Stieltjes constants $\gamma_n(u)$ are the coefficients in the Laurent expansion of the Hurwitz zeta function $\varsigma(s,u)$ about $s=1$

(1.18) $$\varsigma(s,u) = \sum_{n=0}^{\infty}\frac{1}{(n+u)^s} = \frac{1}{s-1} + \sum_{n=0}^{\infty}\frac{(-1)^n}{n!}\gamma_n(u)(s-1)^n$$

and $\gamma_0(u) = -\psi(u)$, where $\psi(u)$ is the digamma function which is the logarithmic derivative of the gamma function $\psi(u) = \frac{d}{du}\log\Gamma(u)$. It is easily seen from the definition of the Hurwitz zeta function that $\varsigma(s,1) = \varsigma(s)$ and accordingly that $\gamma_n(1) = \gamma_n$.

The generalised Euler-Mascheroni constants $\gamma_n$ (or Stieltjes constants) are the coefficients of the Laurent expansion of the Riemann zeta function $\varsigma(s)$ about $s=1$

(1.19) $$\varsigma(s) = \frac{1}{s-1} + \sum_{n=0}^{\infty}\frac{(-1)^n}{n!}\gamma_n(s-1)^n$$

Since $\lim_{s\to 1}\left[\varsigma(s) - \frac{1}{s-1}\right] = \gamma$ it is clear that $\gamma_0 = \gamma$. It may be shown, as in [23, p.4], that

$$(1.20) \qquad \gamma_n = \lim_{N \to \infty} \left[ \sum_{k=1}^{N} \frac{\log^n k}{k} - \frac{\log^{n+1} N}{n+1} \right] = \lim_{N \to \infty} \left[ \sum_{k=1}^{N} \frac{\log^n k}{k} - \int_{1}^{N} \frac{\log^n t}{t} \, dt \right]$$

where we define $\log^0 1 = 1$. The Stieltjes constants may be compared with the Sitaramachandrarao constants $\delta_n$ which are considered later in (2.15).

It was previously shown in [16] that

$$(1.21) \qquad \gamma_n(u) = -\frac{1}{n+1} \sum_{i=0}^{\infty} \frac{1}{i+1} \sum_{j=0}^{i} \binom{i}{j} (-1)^j \log^{n+1}(u+j)$$

## 2. Application of the Lehmer constants $b_n$

Upon equating (1.6) and (1.7.2) we obtain

$$\log \Gamma\left(1 + \frac{s}{2}\right) - \frac{s}{2} \log \pi + \log[(s-1)\varsigma(s)] = -\log 2 - \sum_{k=1}^{\infty} \frac{\sigma_k}{k} s^k$$

and using the Maclaurin expansion [7, p.201]

$$\log \Gamma(1+s) = -\gamma s + \sum_{k=2}^{\infty} (-1)^k \frac{\varsigma(k)}{k} s^k \qquad , -1 < s \leq 1$$

we obtain

$$\log[(s-1)\varsigma(s)] = -\log 2 - \sum_{k=1}^{\infty} \frac{\sigma_k}{k} s^k + \frac{1}{2} \gamma s + \sum_{k=2}^{\infty} (-1)^{k+1} \frac{\varsigma(k)}{k 2^k} s^k + \frac{s}{2} \log \pi$$

This may be expressed as

$$\log[(s-1)\varsigma(s)] = -\log 2 - \sum_{k=1}^{\infty} \frac{\sigma_k}{k} s^k + \sum_{k=1}^{\infty} (-1)^{k+1} \frac{\varsigma(k)}{k 2^k} s^k + \frac{s}{2} \log \pi$$

where, for convenience, we denote $\varsigma(1) = \gamma$.

We then have

$$(2.1) \qquad \log[(s-1)\varsigma(s)] = -\log 2 + \sum_{k=1}^{\infty} \frac{[(-1)^{k+1} 2^{-k} \varsigma(k) - \sigma_k + \frac{1}{2}\delta_{1,k} \log \pi]}{k} s^k$$

and differentiation results in

$$\sum_{k=1}^{\infty}\left[(-1)^{k+1}2^{-k}\varsigma(k)-\sigma_k+\frac{1}{2}\delta_{1,k}\log\pi\right](k-1)(k-1)...(k-m+1)\,s^{k-m}=\frac{d^m}{ds^m}\log[(s-1)\varsigma(s)]$$

Letting $s=0$ gives us

$$\left[(-1)^{m+1}2^{-m}\varsigma(m)-\sigma_m+\frac{1}{2}\delta_{1,m}\log\pi\right](m-1)!=\frac{d^m}{ds^m}\log[(s-1)\varsigma(s)]\Bigg|_{s=0}$$

and we make the finite summation

$$\sum_{m=1}^{n}(-1)^m\binom{n}{m}\left[(-1)^{m+1}2^{-m}\varsigma(m)-\sigma_m+\frac{1}{2}\delta_{1,m}\log\pi\right]$$

$$=\sum_{m=1}^{n}(-1)^m\binom{n}{m}\frac{1}{(m-1)!}\frac{d^m}{ds^m}\log[(s-1)\varsigma(s)]\Bigg|_{s=0}$$

or equivalently

$$-\sum_{m=1}^{n}\binom{n}{m}2^{-m}\varsigma(m)-\sum_{m=1}^{n}(-1)^m\binom{n}{m}\sigma_m-\frac{n}{2}\log\pi=\sum_{m=1}^{n}(-1)^m\binom{n}{m}\frac{1}{(m-1)!}\frac{d^m}{ds^m}\log[(s-1)\varsigma(s)]\Bigg|_{s=0}$$

Then using (1.8) we obtain

$$\lambda_n=\frac{n}{2}\log\pi+\sum_{m=1}^{n}\binom{n}{m}2^{-m}\varsigma(m)+\sum_{m=1}^{n}(-1)^m\binom{n}{m}\frac{1}{(m-1)!}\frac{d^m}{ds^m}\log[(s-1)\varsigma(s)]\Bigg|_{s=0}$$

or equivalently

$$(2.2)\qquad\lambda_n=\frac{n}{2}(\gamma+\log\pi)+\sum_{m=2}^{n}\binom{n}{m}2^{-m}\varsigma(m)+\sum_{m=1}^{n}(-1)^m\binom{n}{m}\frac{1}{(m-1)!}\frac{d^m}{ds^m}\log[(s-1)\varsigma(s)]\Bigg|_{s=0}$$

Prima facie, this formula seems more attractive than (1.6) because we no longer have to deal with the alternating sum $S_1(n)$ defined by (1.16.1).

Lehmer [28] considered the constants $b_n$ defined by

$$(2.3)\qquad\frac{d}{ds}\log[2(s-1)\varsigma(s)]=\frac{\varsigma'(s)}{\varsigma(s)}+\frac{1}{s-1}=\sum_{n=0}^{\infty}b_n s^n\quad,\ |s|<2$$

so that

$$(2.3.1) \qquad \log[2(s-1)\varsigma(s)] = \sum_{n=0}^{\infty} \frac{b_n}{n+1} s^{n+1} = \sum_{n=1}^{\infty} \frac{b_{n-1}}{n} s^n$$

and we note that

$$(2.3.2) \qquad \frac{d^m}{ds^m} \log[2(s-1)\varsigma(s)] \Bigg|_{s=0} = (m-1)! b_{m-1}$$

We then have from (2.2)

$$\lambda_n = \frac{n}{2}(\gamma + \log \pi) + \sum_{m=2}^{n} \binom{n}{m} 2^{-m} \varsigma(m) + \sum_{m=1}^{n} (-1)^m \binom{n}{m} b_{m-1}$$

and since $b_0 = \log 2\pi - 1$ we have

$$\lambda_n = \frac{n}{2}(2 + \gamma - \log 4\pi) + \sum_{m=2}^{n} \binom{n}{m} [2^{-m} \varsigma(m) + (-1)^m b_{m-1}]$$

We know that

$$2 + \gamma - \log 4\pi = 2\lambda_1 \approx 0.046...$$

and therefore we obtain

$$(2.4) \qquad \lambda_n = n\lambda_1 + \sum_{m=2}^{n} \binom{n}{m} [2^{-m} \varsigma(m) + (-1)^m b_{m-1}]$$

We have Lehmer's relation [28] for $m \geq 2$

$$(2.5) \qquad b_{m-1} = (-1)^{m-1} 2^{-m} \varsigma(m) - \sigma_m$$

and noting that

$$\sum_{m=2}^{n} (-1)^m \binom{n}{m} b_{m-1} = -\sum_{m=2}^{n} \binom{n}{m} 2^{-m} \varsigma(m) - \sum_{m=2}^{n} (-1)^m \binom{n}{m} \sigma_m$$

we easily see that

$$\lambda_n = \frac{n}{2}(2 + \gamma - \log 4\pi) - \sum_{m=2}^{n} (-1)^m \binom{n}{m} \sigma_m$$

This may be written as

$$\lambda_n = \frac{n}{2}(2 + \gamma - \log 4\pi) - n\sigma_1 - \sum_{m=1}^{n}(-1)^m \binom{n}{m}\sigma_m$$

which simply leads us back to square one, i.e. equation (1.8)!

Coffey [13] has shown that $b_m$ has strict sign alteration, i.e.

(2.6)        $b_m = (-1)^m \mu_m$ where $\mu_m > 0$

This strict sign alteration was also reported in [17] where we considered the function $L(s) = \log[(s-1)\varsigma(s)]$. We have from (1.12)

$$L^{(1)}(s) = -\sum_{k=0}^{\infty} \eta_k (s-1)^k$$

and hence

$$L^{(1)}(1) = -\eta_0 = \gamma$$

$$L^{(n+1)}(s) = -\sum_{k=0}^{\infty} \eta_k k(k-1)...(k-n+1)(s-1)^{k-n}$$

As we shall see later in Section 3 of this paper, Coffey [13] has also shown that the sequence $(\eta_n)$ has strict sign alteration

$$\eta_n = (-1)^{n+1} \varepsilon_n$$

where $\varepsilon_n$ are positive constants and therefore we have

$$L^{(n+1)}(1) = -n!\eta_n = (-1)^n n!\varepsilon_n$$

We see that

$$L^{(1)}(0) = -\sum_{k=0}^{\infty}(-1)^k \eta_k = \sum_{k=0}^{\infty}\varepsilon_k$$

is positive. In fact we have

$$L^{(1)}(0) = -\sum_{k=0}^{\infty}(-1)^k \eta_k = \log(2\pi) - 1$$

$$L^{(n+1)}(0) = -\sum_{k=0}^{\infty} \eta_k k(k-1)...(k-n+1)(-1)^{k-n} = (-1)^n \sum_{k=0}^{\infty} \varepsilon_k k(k-1)...(k-n+1)$$

and therefore we note that the signs of $L^{(n+1)}(0)$ also strictly alternate.

Coffey [13] has indicated that for large values of $m$ we have

$$\mu_m \approx 2^{-m-1}$$

However, since $\lim_{m \to \infty} \sigma_m = 0$, reference to (2.5) indicates that a better approximation would appear to be

(2.7)     $\mu_{m-1} \approx 2^{-m} \varsigma(m)$

where of course we note that $\lim_{m \to \infty} \varsigma(m) = 1$.

In order to verify the Riemann Hypothesis, equation (2.4) tells us that it would be sufficient to show that

$$\sum_{m=2}^{n} \binom{n}{m} [2^{-m} \varsigma(m) + (-1)^m b_{m-1}] > 0$$

or

$$\sum_{m=2}^{n} \binom{n}{m} [2^{-m} \varsigma(m) - \mu_{m-1}] > 0$$

and indeed it would be sufficient if we could show that $2^{-m} \varsigma(m) > \mu_{m-1}$ for $m \geq 2$. The approximation (2.7) therefore prima facie lends a degree of support to the validity of the Riemann Hypothesis.

However, from (2.5) we have

$$2^{-m} \varsigma(m) - \mu_{m-1} = (-1)^{m+1} \sigma_m$$

and, as reported by Coffey [12, p.16], the initial sign pattern of the $(\sigma_m)$ sequence is simply $--++--\cdots$ with $\sigma_1 > 0$. Therefore the inequality $2^{-m} \varsigma(m) > \mu_{m-1}$ cannot be valid for all $m$. The first 26 values of $\sigma_m$ are reported in Lehmer's paper [28]. The magnitude of $\sigma_{26}$ is given as

$$\sigma_{26} \approx -0.0000\ 0000\ 0000\ 0000\ 0000\ 0000\ 0000\ 01$$

which gives an indication of just how small are the quantities we are dealing with.

Using (2.5) we see that

$$\sum_{m=2}^{\infty} b_{m-1} = \sum_{m=2}^{\infty} \frac{(-1)^{m-1} \varsigma(m)}{2^m} - \sum_{m=2}^{\infty} \sigma_m$$

and it is well known that [13]

$$\sigma_1 = -\sum_{m=1}^{\infty} \sigma_m$$

which gives us

$$2\sigma_1 = -\sum_{m=2}^{\infty} \sigma_m$$

We then have

$$\sum_{m=2}^{\infty} \left[ b_{m-1} - \frac{(-1)^{m-1} \varsigma(m)}{2^m} \right] = 2\sigma_1$$

or

$$\sum_{m=2}^{\infty} (-1)^{m-1} \left[ \mu_{m-1} - \frac{\varsigma(m)}{2^m} \right] = 2\sigma_1$$

which also demonstrates that

$$\lim_{m \to \infty} \left[ \mu_{m-1} - \frac{\varsigma(m)}{2^m} \right] = 0$$

$\square$

We now consider Cauchy's inequality [35, p.84]. If

$$f(s) = \sum_{m=0}^{\infty} a_m s^m \qquad |s| < R$$

and $M(r)$ is the upper bound of $|f(s)|$ on the circle $|s| = r \, (r < R)$ then

$$|a_m| r^m \le M(r)$$

for all values of $m$.

We consider the particular function

$$f(s) = \log[2(s-1)\varsigma(s)] \qquad f(0) = 0$$

We have

$$\left|(s-1)\varsigma(s)\right| = \left|(s-1)\right|\left|\varsigma(s)\right|$$

and for $|s| = r > 1$ we have

$$\left|\varsigma(s)\right| \le \sum_{n=1}^{\infty}\left|\frac{1}{n^s}\right| = \varsigma(r)$$

With $s = re^{i\theta}$ we have

$$\left|(s-1)\right| = \sqrt{r^2 - 2r\cos\theta + 1}$$

so that on the circle $|s| = r$

$$\max\left|(s-1)\right| \le r+1$$

Hence we have

$$\left|f(s)\right| \le \log[2(r+1)\varsigma(r)]$$

and, using $x > \log x$, we obtain

$$\left|f(s)\right| < 2(r+1)\varsigma(r)$$

$$\left|\frac{(-1)^{m-1}\mu_{m-1}}{m}\right| = \frac{\mu_{m-1}}{m} < \frac{2}{r^m}(r+1)\varsigma(r)$$

$$\frac{\mu_{m-1}}{m} < \frac{2}{r^m}(r+1)\varsigma(r) \text{ where } r < 2$$

Unfortunately, this inequality is not sharp enough for our purposes.

$\square$

In this section we show another more concise derivation of (2.4).We have from (1.4.1)

$$\lambda_n = \sum_{m=1}^{n}\binom{n}{m}\frac{1}{(m-1)!}\frac{d^m}{ds^m}\log\xi(s)\Big|_{s=1}$$

Since $\xi(s) = \xi(1-s)$ this is equivalent to

$$= \sum_{m=1}^{n} \binom{n}{m} \frac{1}{(m-1)!} \frac{d^m}{ds^m} \log \xi(1-s) \bigg|_{s=1}$$

Making the substitution $p = 1 - s$ this becomes

$$= \sum_{m=1}^{n} \binom{n}{m} \frac{(-1)^m}{(m-1)!} \frac{d^m}{dp^m} \log \xi(p) \bigg|_{p=0}$$

and making another substitution $p = -s$ we obtain

$$(2.8) \quad \lambda_n = \sum_{m=1}^{n} \binom{n}{m} \frac{1}{(m-1)!} \frac{d^m}{ds^m} \log \xi(-s) \bigg|_{s=0}$$

Letting $s \rightarrow -s$ in (1.6) gives us

$$\log \xi(-s) = \log \Gamma\left(1 - \frac{s}{2}\right) + \frac{s}{2} \log \pi + \log[-(s+1)\varsigma(-s)]$$

Differentiation gives us

$$\frac{d}{ds} \log \Gamma\left(1 - \frac{s}{2}\right) = -\frac{1}{2} \psi\left(1 - \frac{s}{2}\right)$$

and we see that

$$\frac{d^m}{ds^m} \log \Gamma\left(1 - \frac{s}{2}\right) = (-1)^m \frac{1}{2^m} \psi^{(m-1)}\left(1 - \frac{s}{2}\right)$$

Specifically we have

$$\frac{d^m}{ds^m} \log \Gamma\left(1 - \frac{s}{2}\right) \bigg|_{s=0} = (-1)^m \frac{1}{2^m} \psi^{(m-1)}(1) = (-1)^m \frac{1}{2^m} \psi^{(m-1)}(1)$$

We have for $m \geq 2$

$$\psi^{(m-1)}(1) = (-1)^m (m-1)! \varsigma(m)$$

$$\frac{d^m}{ds^m}\log\Gamma\left(1-\frac{s}{2}\right)\bigg|_{s=0} = \frac{1}{2^m}(m-1)!\varsigma(m)$$

$$\frac{d}{ds}\log\Gamma\left(1-\frac{s}{2}\right)\bigg|_{s=0} = \frac{1}{2}\gamma$$

Letting $s \to -s$ in (2.3.1) we see that

$$\log[-2(s+1)\varsigma(-s)] = \sum_{n=1}^{\infty}\frac{(-1)^n b_{n-1}}{n}s^n$$

and we note that

$$\frac{d^m}{ds^m}\log[-2(s+1)\varsigma(-s)]\bigg|_{s=0} = (-1)^m(m-1)!b_{m-1}$$

Dealing separately with the term involving $m=1$ we obtain (2.4) again

$$\lambda_n = \frac{n}{2}\log\pi + \frac{n}{2}\gamma - n(\log 2\pi - 1) + \sum_{m=2}^{n}\binom{n}{m}\left[\frac{\varsigma(m)}{2^m} + (-1)^m b_{m-1}\right]$$

☐

From (2.3) we have

$$\frac{d}{ds}[2(s-1)\varsigma(s)] = 2(s-1)\varsigma(s)\sum_{n=0}^{\infty}b_n s^n$$

and hence referring to (10.10) we have

(2.9) $$\frac{d^m}{ds^m}[2(s-1)\varsigma(s)]\bigg|_{s=0} = Y_m(0!b_0, 1!b_1, ..., (m-1)!b_{m-1})$$

in terms of the (exponential) complete Bell polynomials.

We have the Maclaurin expansion

(2.10) $$\varsigma(s) = \sum_{n=0}^{\infty}\varsigma^{(n)}(0)\frac{s^n}{n!}$$

and we see that

$$\sum_{n=0}^{\infty} \varsigma^{(n-1)}(0) \frac{ns^n}{n!} = \sum_{n=1}^{\infty} \varsigma^{(n-1)}(0) \frac{ns^n}{n!} = s \sum_{n=1}^{\infty} \varsigma^{(n-1)}(0) \frac{s^{n-1}}{(n-1)!}$$

and in turn we have

$$s \sum_{n=1}^{\infty} \varsigma^{(n-1)}(0) \frac{s^{n-1}}{(n-1)!} = s \sum_{m=0}^{\infty} \varsigma^{(m)}(0) \frac{s^m}{m!} = s\varsigma(s)$$

Hence we have

$$(s-1)\varsigma(s) = \sum_{n=0}^{\infty} \left[ n\varsigma^{(n-1)}(0) - \varsigma^{(n)}(0) \right] \frac{s^n}{n!}$$

and differentiation shows that

$$\frac{d^m}{ds^m} [2(s-1)\varsigma(s)] \bigg|_{s=0} = 2 \left[ m\varsigma^{(m-1)}(0) - \varsigma^{(m)}(0) \right]$$

Therefore we obtain

(2.11) $$2 \left[ m\varsigma^{(m-1)}(0) - \varsigma^{(m)}(0) \right] = Y_m(0!b_0, 1!b_1, ..., (m-1)!b_{m-1})$$

and for example with $m=1$ we have

$$2 \left[ \varsigma(0) - \varsigma^{(1)}(0) \right] = Y_1(0!b_0) = b_0$$

so that

$$b_0 = \log 2\pi - 1$$

With $m=2$ we obtain

$$2 \left[ 2\varsigma^{(1)}(0) - \varsigma^{(2)}(0) \right] = Y_2(0!b_0, 1!b_1) = b_0^2 + b_1$$

so that

(2.12) $$b_1 = 2 \left[ 2\varsigma^{(1)}(0) - \varsigma^{(2)}(0) \right] - 4 \left[ \varsigma(0) - \varsigma^{(1)}(0) \right]^2$$

Since $b_1$ is negative we see that

$$2 \left[ \varsigma(0) - \varsigma^{(1)}(0) \right]^2 > \left[ 2\varsigma^{(1)}(0) - \varsigma^{(2)}(0) \right]$$

With $s \to 1-s$ in (2.10) we have

$$\varsigma(1-s) = \sum_{n=0}^{\infty} \frac{\varsigma^{(n)}(0)}{n!}(1-s)^n$$

and combining this with the geometric series

$$\frac{1}{s} = \sum_{n=0}^{\infty}(1-s)^n$$

we obtain

$$(2.13) \quad \varsigma(1-s) + \frac{1}{s} = \sum_{n=0}^{\infty}\left[\frac{\varsigma^{(n)}(0)}{n!}+1\right](1-s)^n$$

With $1-s = p$ this becomes

$$\varsigma(p) + \frac{1}{1-p} = \sum_{n=0}^{\infty}\left[\frac{\varsigma^{(n)}(0)}{n!}+1\right]p^n$$

and we see the connection with the $\delta_n$ constants considered by Sitaramachandrarao [32] in 1986

$$(2.14) \quad \varsigma(s) + \frac{1}{1-s} = \sum_{n=0}^{\infty}\frac{(-1)^n \delta_n}{n!} s^n$$

where

$$(2.15) \quad \delta_n = \lim_{m\to\infty}\left[\sum_{k=1}^{m}\log^n k - \int_{1}^{m}\log^n x\,dx - \frac{1}{2}\log^n m\right]$$

$$= (-1)^n\left[\varsigma^{(n)}(0) + n!\right]$$

A derivation of (2.15) follows. Using the Euler-Maclaurin summation formula [3], Hardy [22, p.333] showed that the Riemann zeta function could be expressed as follows

$$(2.16) \quad \varsigma(s) = \lim_{n\to\infty}\left[\sum_{k=1}^{n}\frac{1}{k^s} - \frac{n^{1-s}}{1-s} - \frac{1}{2}n^{-s}\right] \qquad \mathrm{Re}(s) > -1$$

It may immediately be seen that this identity is trivially satisfied for $\mathrm{Re}(s) > 1$ because

$$\varsigma(s) = \lim_{n \to \infty} \left[ \sum_{k=1}^{n} \frac{1}{k^s} - \frac{n^{1-s}}{1-s} - \frac{1}{2} n^{-s} \right] = \lim_{n \to \infty} \sum_{k=1}^{n} \frac{1}{k^s} + \lim_{n \to \infty} \left[ -\frac{n^{1-s}}{1-s} - \frac{1}{2} n^{-s} \right]$$

and the latter limit is clearly equal to zero.

Differentiating (2.16) results in for $\operatorname{Re}(s) > -1$

(2.17) $$\varsigma'(s) = \lim_{n \to \infty} \left[ -\sum_{k=1}^{n} \frac{\log k}{k^s} + \frac{n^{1-s}(1-s)\log n - n^{1-s}}{(1-s)^2} + \frac{1}{2} n^{-s} \log n \right]$$

and with $s = 0$ we obtain

$$\varsigma'(0) = \lim_{n \to \infty} \left[ -\sum_{k=1}^{n} \log k + \left( n + \frac{1}{2} \right) \log n - n \right]$$

Hence, using the Stirling approximation for $n!$ we see that $\varsigma'(0) = -\frac{1}{2} \log(2\pi)$.

With regard to (2.16) we could determine $\varsigma''(0)$

$$\varsigma''(s) = \lim_{n \to \infty} \left[ \sum_{k=1}^{n} \frac{\log^2 k}{k^s} + \frac{(1-s)^2[-n^{1-s}(1-s)\log^2 n] + 2[n^{1-s}(1-s)\log n - n^{1-s}]}{(1-s)^4} - \frac{1}{2} n^{-s} \log^2 n \right]$$

so that

(2.18) $$\varsigma''(0) = \lim_{n \to \infty} \left[ \sum_{k=1}^{n} \log^2 k - n \log^2 n + 2n \log n - 2n - \frac{1}{2} \log^2 n \right]$$

and compare the result with

$$\varsigma''(0) = \gamma_1 + \frac{1}{2} \gamma^2 - \frac{1}{24} \pi^2 - \frac{1}{2} \log^2(2\pi)$$

The equation (2.18) is contained in Ramanujan's Notebooks [5, Part I, p.203].

In order to simplify the calculations, we write

$$\varsigma(s) = \lim_{n \to \infty} \left[ \sum_{k=1}^{n} \frac{1}{k^s} - \frac{n^{1-s}}{1-s} - \frac{1}{2} n^{-s} \right] = \lim_{n \to \infty} \left[ \sum_{k=1}^{n} \frac{1}{k^s} - \frac{n^{1-s}-1}{1-s} - \frac{1}{1-s} - \frac{1}{2} n^{-s} \right]$$

and differentiation gives us

$$\varsigma^{(n)}(s) = \lim_{m \to \infty}\left[(-1)^n \sum_{k=1}^{m} \frac{\log^n k}{k^s} - f^{(n)}(s) - \frac{n!}{(1-s)^{n+1}} - \frac{1}{2}(-1)^n m^{-s} \log^n m\right]$$

where we have denoted $f(s)$ as

$$f(s) = \frac{m^{1-s} - 1}{s-1}$$

We can represent $f(s)$ by the following integral

$$f(s) = \frac{m^{1-s} - 1}{s-1} = -\int_{1}^{m} x^{-s} dx$$

so that

$$f^{(n)}(s) = -(-1)^n \int_{1}^{m} x^{-s} \log^n x\, dx$$

and thus

$$f^{(n)}(0) = -(-1)^n \int_{1}^{m} \log^n x\, dx$$

Therefore, with $s = 0$ we obtain

(2.19)    $$(-1)^n \left[\varsigma^{(n)}(0) + n!\right] = \lim_{m \to \infty}\left[\sum_{k=1}^{m} \log^n k - \int_{1}^{m} \log^n x\, dx - \frac{1}{2}\log^n m\right]$$

As pointed out by Apostol [2], the power series expansion (2.13) converges for $s = 0$ and therefore we have

(2.20)    $$\lim_{n \to \infty}\left[\frac{\varsigma^{(n)}(0)}{n!} + 1\right] = 0$$

Apostol [2] has calculated the first 18 values of $\frac{\varsigma^{(n)}(0)}{n!}$ and we note that they exhibit some small oscillations around the value $-1$.

We also see that

$$\lim_{n \to \infty}\left[\frac{\varsigma^{(n-1)}(0)}{(n-1)!} + 1\right] - \lim_{n \to \infty}\left[\frac{\varsigma^{(n)}(0)}{n!} + 1\right] = 0$$

and thus

$$\lim_{n \to \infty} \frac{1}{n!} \left[ n\varsigma^{(n-1)}(0) - \varsigma^{(n)}(0) \right] = 0$$

We note that

$$\lim_{s \to 1} \left[ \varsigma(s) - \frac{1}{s-1} \right] = \sum_{n=0}^{\infty} \left[ \frac{\varsigma^{(n)}(0)}{n!} + 1 \right] = \gamma$$

We have the Laurent expansion of the zeta function $\varsigma(s)$ about $s = 1$

$$\varsigma(s) = \frac{1}{s-1} + \sum_{n=0}^{\infty} \frac{(-1)^n}{n!} \gamma_n (s-1)^n$$

where $\gamma_n$ are known as the Stieltjes constants. We have $\gamma_0 = -\psi(1) = \gamma$.

We see that

$$\gamma = \sum_{n=0}^{\infty} \left[ \frac{\varsigma^{(n)}(0)}{n!} + 1 \right] = \frac{1}{2} + \sum_{n=1}^{\infty} \left[ \frac{\varsigma^{(n)}(0)}{n!} + 1 \right]$$

and more generally we obtain by differentiation for $p \geq 1$

$$(-1)^p \gamma_p = \sum_{n=p}^{\infty} \left[ \frac{\varsigma^{(n)}(0)}{n!} + 1 \right] n(n-1)...(n-p+1)$$

which may be written as

$$(2.21) \qquad \frac{(-1)^p}{p!} \gamma_p = \sum_{n=p}^{\infty} \binom{n}{p} \left[ \frac{\varsigma^{(n)}(0)}{n!} + 1 \right]$$

We showed in [16] that

$$(2.22) \qquad (-1)^r \varsigma^{(r)}(s,u) = \frac{r!}{(s-1)^{r+1}} + \sum_{p=0}^{\infty} \frac{(-1)^p}{p!} (s-1)^p \gamma_{p+r}(u)$$

In 1995, Choudhury [10] mentioned that Ramanujan determined that for $r \geq 1$ and $\mathrm{Re}(s) > 1$

$$(-1)^r \varsigma^{(r)}(s) = \sum_{p=1}^{\infty} \frac{\log^r p}{p^s} = \frac{r!}{(s-1)^{r+1}} + \sum_{p=0}^{\infty} \frac{(-1)^p}{p!} (s-1)^p \gamma_{p+r}$$

(see also Ramanujan's Notebooks [21, Part I, p.224]). This is a particular case of (2.22) where $u = 1$ and it should be noted that (2.22) is valid for all $s \neq 1$.

We have for $s = 0$

$$\varsigma^{(n)}(0) = -n! + (-1)^n \sum_{p=0}^{\infty} \frac{\gamma_{p+n}}{p!}$$

This shows that the Sitaramachandrarao constants $\delta_n$ may be represented by

$$(2.23) \qquad \delta_n = (-1)^n \left[ \varsigma^{(n)}(0) + n! \right] = \sum_{p=0}^{\infty} \frac{\gamma_{p+n}}{p!}$$

We see that

$$\varsigma^{(n-1)}(0) = -(n-1)! + (-1)^{n-1} \sum_{p=0}^{\infty} \frac{\gamma_{p+n-1}}{p!}$$

$$-\varsigma^{(n)}(0) + n\varsigma^{(n-1)}(0) = (-1)^{n-1} n \sum_{p=0}^{\infty} \frac{\gamma_{p+n-1}}{p!} - (-1)^n \sum_{p=0}^{\infty} \frac{\gamma_{p+n}}{p!}$$

$$-\varsigma^{(n)}(0) + n\varsigma^{(n-1)}(0) = (-1)^{n+1} \sum_{p=0}^{\infty} \frac{n\gamma_{p+n-1} + \gamma_{p+n}}{p!}$$

and reference to (2.11) shows that

$$(2.24) \qquad Y_n(0!b_0, 1!b_1, ..., (n-1)!b_{n-1}) = 2(-1)^{n+1} \sum_{p=0}^{\infty} \frac{n\gamma_{p+n-1} + \gamma_{p+n}}{p!}$$

$\square$

We may write (2.14) as

$$(s-1)\varsigma(s) = 1 + \sum_{n=0}^{\infty} \frac{(-1)^n \delta_n}{n!} \left[ s^{n+1} - s^n \right]$$

We have

$$\frac{d^m}{ds^m}(s-1)\varsigma(s) = m\varsigma^{(m-1)}(s) + (s-1)\varsigma^{(m)}(s)$$

$$= \sum_{n=0}^{\infty} \frac{(-1)^n \delta_n}{n!} \left[ (n+1)...(n-m)s^{n+1-m} - n...(n-m+1)s^{n-m} \right]$$

and evaluation at $s = 0$ results in

$$m\varsigma^{(m-1)}(0) - \varsigma^{(m)}(0) = (-1)^{m-1}[m\delta_{m-1} + \delta_m]$$

We recall Hasse's formula for the zeta function that for $\operatorname{Re}(s) \neq 1$ (see for example [15])

$$(s-1)\varsigma(s) = \sum_{n=0}^{\infty} \frac{1}{n+1} \sum_{k=0}^{n} \binom{n}{k} \frac{(-1)^k}{(1+k)^{s-1}}$$

and we see that

$$\frac{d^m}{ds^m}(s-1)\varsigma(s) = (-1)^m \sum_{n=0}^{\infty} \frac{1}{n+1} \sum_{k=0}^{n} \binom{n}{k} \frac{(-1)^k \log^m(1+k)}{(1+k)^{s-1}}$$

Evaluation at $s = 0$ gives us

$$m\varsigma^{(m-1)}(0) - \varsigma^{(m)}(0) = (-1)^m \sum_{n=0}^{\infty} \frac{1}{n+1} \sum_{k=0}^{n} \binom{n}{k} (-1)^k (1+k) \log^m(1+k)$$

$$= (-1)^m \sum_{n=0}^{\infty} \frac{1}{n+1} \sum_{k=0}^{n} \binom{n}{k} (-1)^k \log^m(1+k) + (-1)^m \sum_{n=0}^{\infty} \frac{1}{n+1} \sum_{k=0}^{n} \binom{n}{k} (-1)^k k \log^m(1+k)$$

We see from (1.21) that

$$\gamma_{m-1} = -\frac{1}{m} \sum_{n=0}^{\infty} \frac{1}{n+1} \sum_{k=0}^{n} \binom{n}{k} (-1)^k \log^m(1+k)$$

and hence we have for $m \geq 1$

$$\sum_{n=0}^{\infty} \frac{1}{n+1} \sum_{k=0}^{n} \binom{n}{k} (-1)^k k \log^m(1+k) = (-1)^m \left[ m\varsigma^{(m-1)}(0) - \varsigma^{(m)}(0) \right] + \frac{\gamma_{m-1}}{m}$$

□

Adamchik [1] noted that the Hermite integral for the Hurwitz zeta function may be derived from the Abel-Plana summation formula [33, p.90]

$$(2.25) \qquad \sum_{k=0}^{\infty} f(k) = \frac{1}{2} f(0) + \int_{0}^{\infty} f(x)\,dx + i \int_{0}^{\infty} \frac{f(ix) - f(-ix)}{e^{2\pi x} - 1}\,dx$$

which applies to functions which are analytic in the right-hand plane and satisfy the convergence condition

$$\lim_{y \to \infty} e^{-2\pi y} \left| f(x+iy) \right| = 0$$

uniformly on any finite interval of $x$. Derivations of the Abel-Plana summation formula are given in [36, p.145] and [37, p.108].

Letting $f(k) = (k+u)^{-s}$ we obtain

$$(2.26) \qquad \varsigma(s,u) = \sum_{k=0}^{\infty} \frac{1}{(k+u)^s} = \frac{u^{-s}}{2} + \frac{u^{1-s}}{s-1} + i\int_0^{\infty} \frac{(u+ix)^{-s} - (u-ix)^{-s}}{e^{2\pi x} - 1} dx$$

Then, noting that

$$(u+ix)^{-s} - (u-ix)^{-s} = (re^{i\theta})^{-s} - (re^{-i\theta})^{-s}$$

$$= r^{-s}[e^{-is\theta} - e^{is\theta}]$$

$$= \frac{2}{i(u^2+x^2)^{s/2}} \sin(s\tan^{-1}(x/u))$$

we may write (2.26) as Hermite's integral for $\varsigma(s,u)$

$$(2.27) \qquad \varsigma(s,u) = \frac{u^{-s}}{2} + \frac{u^{1-s}}{s-1} + 2\int_0^{\infty} \frac{\sin(s\tan^{-1}(x/u))}{(u^2+x^2)^{s/2}(e^{2\pi x}-1)} dx$$

We now take one step back and differentiate the intermediate equation (2.26) in the case where $u=1$ with respect to $s$ to obtain

$$(2.28) \qquad \varsigma^{(n)}(s) = \frac{(-1)^n n!}{(s-1)^{n+1}} + i(-1)^n \int_0^{\infty} \frac{(1+ix)^{-s}\log^n(1+ix) - (1-ix)^{-s}\log^n(1-ix)}{e^{2\pi x} - 1} dx$$

It may be noted that

$$i(-1)^n \int_0^{\infty} \frac{(1+ix)^{-s}\log^n(1+ix) - (1-ix)^{-s}\log^n(1-ix)}{e^{2\pi x} - 1} dx = 2(-1)^{n+1} \operatorname{Im} \int_0^{\infty} \frac{(1+ix)^{-s}\log^n(1+ix)}{e^{2\pi x} - 1} dx$$

and with $s=0$ we obtain

$$\varsigma^{(n)}(0) = -n! + 2(-1)^{n+1} \operatorname{Im} \int_0^\infty \frac{\log^n(1+ix)}{e^{2\pi x}-1} dx$$

$$\frac{\varsigma^{(n)}(0)}{n!} + 1 = 2\frac{(-1)^{n+1}}{n!} \operatorname{Im} \int_0^\infty \frac{\log^n(1+ix)}{e^{2\pi x}-1} dx$$

In particular we have

$$\varsigma^{(1)}(0) + 1 = 2\operatorname{Im} \int_0^\infty \frac{\log(1+ix)}{e^{2\pi x}-1} dx$$

$$= 2\int_0^\infty \frac{\tan^{-1}x}{e^{2\pi x}-1} dx$$

which shows that

$$\varsigma^{(1)}(0) + 1 > 0$$

We also have

$$\frac{\varsigma^{(2)}(0)}{2!} + 1 = -\operatorname{Im} \int_0^\infty \frac{\log^2(1+ix)}{e^{2\pi x}-1} dx$$

$$= -\int_0^\infty \frac{\log\left(1+x^2\right)\tan^{-1}x}{e^{2\pi x}-1} dx$$

which shows that

$$0 > \frac{\varsigma^{(2)}(0)}{2!} + 1$$

The above inequalities concur with the numerical values in [2]. However, it is not immediately obvious how this procedure may be extended to higher values of $n$.

## 3. Application of the eta constants $\eta_n$

We have from (1.1)

$$\frac{d^m}{ds^m}\log[(s-1)\varsigma(s)] = -\sum_{k=1}^\infty \frac{\eta_{k-1}}{k}k(k-1)(k-2)...(k-m+1)(s-1)^{k-m}$$

so that

$$\frac{d^m}{ds^m}\log[(s-1)\varsigma(s)]\Bigg|_{s=0} = -\sum_{k=1}^{\infty}\frac{\eta_{k-1}}{k}k(k-2)(k-1)...(k-m+1)(-1)^{k-m}$$

$$= (-1)^{m+1}\sum_{k=1}^{\infty}\frac{\eta_{k-1}}{k}(-1)^k k(k-1)(k-2)...(k-m+1)$$

Therefore we have

(3.1)     $$\sum_{m=1}^{n}(-1)^m\binom{n}{m}\frac{1}{(m-1)!}\frac{d^m}{ds^m}\log[(s-1)\varsigma(s)]\Bigg|_{s=0}$$

$$= -\sum_{m=1}^{n}\binom{n}{m}\frac{1}{(m-1)!}\sum_{k=1}^{\infty}(-1)^k\eta_{k-1}(k-1)(k-2)...(k-m+1)$$

Coffey [13] has shown that the sequence $(\eta_k)$ has strict sign alteration

(3.2)     $$\eta_k = (-1)^{k+1}\varepsilon_k$$

where $\varepsilon_k$ are positive constants. This sign alteration was also noted by Israilov. With reference to (2.2) we may therefore write

(3.3)     $$\sum_{m=1}^{n}(-1)^m\binom{n}{m}\frac{1}{(m-1)!}\frac{d^m}{ds^m}\log[(s-1)\varsigma(s)]\Bigg|_{s=0}$$

$$= -\sum_{m=1}^{n}\binom{n}{m}\frac{1}{(m-1)!}\sum_{k=1}^{\infty}\varepsilon_{k-1}(k-1)(k-2)...(k-m+1)$$

It is certainly a cruel twist of fate that a solitary negative sign appears at this stage of the analysis because reference to (2.2) shows that a plus sign here instead would have immediately resulted in the Riemann Hypothesis being verified!

Let us now consider a specific example.

$$\lambda_1 = \frac{1}{2}(\gamma + \log\pi) - \frac{d}{ds}\log[(s-1)\varsigma(s)]\Bigg|_{s=0}$$

$$= \frac{1}{2}(\gamma + \log\pi) - \frac{[-\varsigma'(0)+\varsigma(0)]}{-\varsigma(0)}$$

$$= \frac{1}{2}(\gamma + \log \pi) - 2[-\varsigma'(0) + \varsigma(0)]$$

Hence we have

(3.4)     $\lambda_1 = \frac{1}{2}\gamma - \frac{1}{2}\log \pi + 1 - \log 2$

where we have used the well-known values

(3.5)     $\varsigma(0) = -\frac{1}{2}$     $\varsigma'(0) = -\frac{1}{2}\log(2\pi)$

□

We may write

$$\log[(s-1)\varsigma(s)] = \log[1 + (s-1)\varsigma(s) - 1]$$

and applying the logarithmic expansion we obtain

$$\log[(s-1)\varsigma(s)] = -\sum_{k=1}^{\infty} \frac{(-1)^k}{k}[(s-1)\varsigma(s) - 1]^k$$

Differentiation of (2.1) gives us

$$\sum_{k=1}^{\infty}\left[(-1)^{k+1}2^{-k}\varsigma(k) - \sigma_k + \frac{1}{2}\delta_{1,k}\log \pi\right]s^{k-1} = -\sum_{k=1}^{\infty}(-1)^k[(s-1)\varsigma(s) - 1]^{k-1}[(s-1)\varsigma'(s) + \varsigma(s)]$$

and with $s = 0$ we obtain

$$\frac{1}{2}\gamma - \sigma_1 + \frac{1}{2}\log \pi = [-\varsigma'(0) + \varsigma(0)]\sum_{k=1}^{\infty}\frac{1}{2^{k-1}}$$

which simplifies to

$$\sigma_1 = \frac{1}{2}\gamma + \frac{1}{2}\log \pi - 2[-\varsigma'(0) + \varsigma(0)]$$

or equivalently

$$\sigma_1 = \frac{1}{2}\gamma - \frac{1}{2}\log \pi + 1 - \log 2$$

A further differentiation results in

$$\sum_{k=1}^{\infty}\left[(-1)^{k+1}2^{-k}\varsigma(k)-\sigma_k+\frac{1}{2}\delta_{1,k}\log\pi\right](k-1)s^{k-2}$$

$$=-\sum_{k=1}^{\infty}(-1)^k(k-1)[(s-1)\varsigma(s)-1]^{k-2}[(s-1)\varsigma'(s)+\varsigma(s)]^2$$

$$-\sum_{k=1}^{\infty}(-1)^k[(s-1)\varsigma(s)-1]^{k-1}[(s-1)\varsigma''(s)+2\varsigma'(s)]$$

and with $s=0$ we obtain

$$-2^{-2}\varsigma(2)-\sigma_2=[-\varsigma'(0)+\varsigma(0)]^2\sum_{k=1}^{\infty}\frac{k-1}{2^{k-2}}+[-\varsigma''(0)+2\varsigma'(0)]\sum_{k=1}^{\infty}\frac{1}{2^{k-1}}$$

Since $\dfrac{1}{1-x}=\displaystyle\sum_{k=1}^{\infty}x^{k-1}$ we have by differentiation

$$\frac{1}{(1-x)^2}=\sum_{k=1}^{\infty}(k-1)x^{k-2}$$

so that $\displaystyle\sum_{k=1}^{\infty}\frac{k-1}{2^{k-2}}=4$ . We then see that

$$-2^{-2}\varsigma(2)-\sigma_2=4[-\varsigma'(0)+\varsigma(0)]^2+2[-\varsigma''(0)+2\varsigma'(0)]$$

Using $\lambda_2=2\sigma_1-\sigma_2$ we obtain

(3.6)     $\lambda_2=\gamma+\log\pi+2^{-2}\varsigma(2)-4[\varsigma(0)-\varsigma'(0)]+4[\varsigma(0)-\varsigma'(0)]^2+2[2\varsigma'(0)-\varsigma''(0)]$

As previously derived by Ramanujan [5] and Apostol [2] we have

(3.7)     $\varsigma''(0)=\gamma_1+\dfrac{1}{2}\gamma^2-\dfrac{1}{4}\varsigma(2)-\dfrac{1}{2}\log^2(2\pi)$

and therefore, by a rather circuitous route, we have obtained the second Li/Keiper constant

(3.8)     $\lambda_2=\dfrac{3}{4}\varsigma(2)+1+\gamma-\gamma^2-2\log 2-\log\pi-2\gamma_1$

□

We recall (1.15)

$$\eta_n = (-1)^{n+1} \left[ \sigma_{n+1} + \left( 1 - \frac{1}{2^{n+1}} \right) \varsigma(n+1) - 1 \right]$$

and using (3.2) $\eta_n = (-1)^{n+1} \varepsilon_n$ we see that

$$\varepsilon_n = \sigma_{n+1} + \left( 1 - \frac{1}{2^{n+1}} \right) \varsigma(n+1) - 1$$

and hence we have

(3.9) $\qquad \sigma_{n+1} > 1 - \left( 1 - \frac{1}{2^{n+1}} \right) \varsigma(n+1)$

We recall Lehmer's relation (2.5) for $n \geq 1$

$$\sigma_{n+1} = (-1)^n 2^{-n-1} \varsigma(n+1) - b_n$$

and we deduce that for $n \geq 1$

(3.10) $\qquad \varsigma(n+1) - 1 - [1 + (-1)^{n+1}] \frac{\varsigma(n+1)}{2^{n+1}} > b_n$

With $n \to 2n$ this inequality becomes

$$\varsigma(2n+1) - 1 > b_{2n}$$

With $n \to 2n - 1$ we obtain

$$\varsigma(2n) \left[ 1 - \frac{1}{2^{2n-1}} \right] - 1 > b_{2n-1}$$

## 4. A conjecture for the Riemann Hypothesis

Coffey [12] noted that the radius of convergence of (1.11)

$$\frac{d}{ds} \log[(s-1)\varsigma(s)] = \frac{\varsigma'(s)}{\varsigma(s)} + \frac{1}{s-1} = -\sum_{k=1}^{\infty} \eta_{k-1}(s-1)^{k-1}$$

is 3 because the first singularity encountered is the trivial zero of $\varsigma(s)$ at $s = -2$. Furthermore he remarked that $|\eta_k|$ cannot increase faster than $3^{-k}$ for sufficiently large $k$ and that $\eta_k \simeq -\gamma(-1/3)^k$ is a very good approximation.

The first 60 values of $\eta_k$ are set out in [12, p.29] and, based on this limited data, it appears that

(4.1) $$\eta_k \simeq (-1/3)^{k+1}$$

is a better approximation for $k \geq 1$. For example, we have

$$\eta_{10} \simeq -5.66605 \times 10^{-6} \qquad (-1/3)^{11} \simeq -5.64503 \times 10^{-6}$$
$$\eta_{20} \simeq -9.56012 \times 10^{-11} \qquad (-1/3)^{21} \simeq -9.55990 \times 10^{-11}$$
$$\eta_{30} \simeq -1.61898 \times 10^{-15} \qquad (-1/3)^{31} \simeq -1.61897 \times 10^{-15}$$
$$\eta_{40} \simeq -2.74176 \times 10^{-20} \qquad (-1/3)^{41} \simeq -2.74175 \times 10^{-20}$$
$$\eta_{50} \simeq -4.64318 \times 10^{-25} \qquad (-1/3)^{51} \simeq -4.64318 \times 10^{-25}$$
$$\eta_{60} \simeq -7.86327 \times 10^{-30} \qquad (-1/3)^{61} \simeq -7.86327 \times 10^{-30}$$

It is possible that $\eta_k \simeq (-1/3)^{k+1} \varsigma(k+1)$ may be an even better approximation, but this aspect has not been pursued any further in this paper.

Referring to (3.1) we consider the summation

$$S = \sum_{k=1}^{\infty} \varepsilon_{k-1}(k-1)(k-2)...(k-m+1)$$

and let us conjecture that for all $k \geq 1$

(4.2) $$|\eta_k| \leq \alpha 3^{-k-1}$$

where $\alpha > 0$, with the result that

$$S \leq \alpha \sum_{k=1}^{\infty} \frac{1}{3^k}(k-1)(k-2)...(k-m+1)$$

Since $\dfrac{1}{1-x} = \sum_{k=1}^{\infty} x^{k-1}$ we have by differentiation

$$\frac{1}{(1-x)^2} = \sum_{k=1}^{\infty} (k-1)x^{k-2}$$

Higher derivatives result in

$$\frac{(m-1)!}{(1-x)^m} = \sum_{k=1}^{\infty} (k-1)(k-2)...(k-m+1)x^{k-m}$$

which may be written as

$$\sum_{k=1}^{\infty} (k-1)(k-2)...(k-m+1)x^k = (m-1)!\left(\frac{x}{1-x}\right)^m$$

With $x = 1/3$ we obtain

$$\sum_{k=1}^{\infty} \frac{1}{3^k}(k-1)(k-2)...(k-m+1) = (m-1)!\frac{1}{2^m}$$

and we therefore have

$$S \leq \alpha(m-1)!\frac{1}{2^m}$$

Hence we see that

$$\sum_{m=1}^{n}\binom{n}{m}\frac{1}{(m-1)!}\sum_{k=1}^{\infty}\varepsilon_{k-1}(k-1)(k-2)...(k-m+1) \leq \alpha\sum_{m=1}^{n}\binom{n}{m}\frac{1}{2^m}$$

Substituting

$$\sum_{m=1}^{n}\binom{n}{m}\frac{1}{2^m} = \sum_{m=0}^{n}\binom{n}{m}\frac{1}{2^m} - 1 = \left(\frac{3}{2}\right)^n - 1$$

gives us

$$\sum_{m=1}^{n}\binom{n}{m}\frac{1}{(m-1)!}\sum_{k=1}^{\infty}\varepsilon_{k-1}(k-1)(k-2)...(k-m+1) \leq \alpha\left[\left(\frac{3}{2}\right)^n - 1\right]$$

We recall (3.1)

$$\lambda_n = \frac{n}{2}(\gamma + \log\pi) + \sum_{m=2}^{n}\binom{n}{m}2^{-m}\varsigma(m) - \sum_{m=1}^{n}\binom{n}{m}\frac{1}{(m-1)!}\sum_{k=1}^{\infty}\varepsilon_{k-1}(k-1)(k-2)...(k-m+1)$$

and we therefore have

$$\lambda_n \geq \frac{n}{2}(\gamma + \log \pi) + \sum_{m=2}^{n} \binom{n}{m} 2^{-m} \varsigma(m) - \alpha \left(\frac{3}{2}\right)^n + \alpha$$

Employing the inequality $\varsigma(m) > 1$ we see that

$$\sum_{m=2}^{n} \binom{n}{m} 2^{-m} \varsigma(m) > \sum_{m=2}^{n} \binom{n}{m} 2^{-m}$$

and using

$$\sum_{m=2}^{n} \binom{n}{m} \frac{1}{2^m} = \sum_{m=0}^{n} \binom{n}{m} \frac{1}{2^m} - 1 - \frac{1}{2} n = \left(\frac{3}{2}\right)^n - 1 - \frac{1}{2} n$$

we obtain

(4.3) $\qquad \lambda_n > \frac{n}{2}(\gamma + \log \pi - 1) + (1-\alpha)\left[\left(\frac{3}{2}\right)^n - 1\right]$

where we note that $\gamma + \log \pi - 1 = 0.72...$ is positive. Therefore the conjecture that $|\eta_k| \leq \alpha 3^{-k-1}$, combined with $1 > \alpha > 0$, immediately implies that $\lambda_n > 0$.

However, it should be noted that the factor $\left(\frac{3}{2}\right)^n$ increases very rapidly; for example $\left(\frac{3}{2}\right)^{100} \approx 10^{17}$ whereas we note from [12] that $\lambda_{100} = 118.6038...$. It is therefore clear that our initial conjecture that $|\eta_k| \leq \alpha 3^{-k-1}$ for **all** $k$ can only be valid provided $(1-\alpha)$ is sufficiently close to zero. It is however clear that this condition is not valid because $\eta_1 \simeq 0.187546 > 1/3^2$.

Further work is required to determine if the conjecture may be modified so that $|\eta_k| \leq \alpha 3^{-k-1}$ for all $k > N$ where $N$ is a fixed integer.

The limited data in [12], also suggests that the terms $|\eta_k|$ form a monotonic decreasing sequence.

## 5. Using the positivity of $\xi^{(n)}(1)$

Using the Jacobi theta function, Coffey [11] showed in 2003 that the even derivatives of $\xi(s)$ are positive for all real values of $s$, while the odd derivatives are positive for $s \geq \frac{1}{2}$ and negative for $s < \frac{1}{2}$. In particular, $\xi^{(n)}(1)$ is positive for $n \geq 1$. This was also proved by Freitas [20] in a different manner in 2005. It may also be noted from the functional equation (1.3) that

(5.1) $$\xi^{(n)}(0) = (-1)^n \xi^{(n)}(1)$$

We write

$$\log[2\xi(s)] = \log\left(1 + [2\xi(s) - 1]\right)$$

and applying the Maclaurin expansion we have

$$\log[2\xi(s)] = -\sum_{k=1}^{\infty} \frac{(-1)^k}{k} \left[2\xi(s) - 1\right]^k$$

Since

$$\frac{d^n}{ds^n} s^{n-1} \log[2\xi(s)] = \frac{d^n}{ds^n} s^{n-1} \log \xi(s)$$

we note from (1.4) that

$$\lambda_n = \frac{1}{(n-1)!} \frac{d^n}{ds^n} [s^{n-1} \log 2\xi(s)] \Bigg|_{s=1}$$

We have as before

$$\frac{d^n}{ds^n} [s^{n-1} \log 2\xi(s)] = \sum_{m=1}^{n} \binom{n}{m} \frac{(n-1)!}{(m-1)!} \frac{d^m}{ds^m} [\log 2\xi(s)]$$

and therefore we get

$$\lambda_n = -\sum_{m=1}^{n} \binom{n}{m} \frac{1}{(m-1)!} \frac{d^m}{ds^m} \sum_{k=1}^{\infty} \frac{(-1)^k}{k} \left[2\xi(s) - 1\right]^k \Bigg|_{s=1}$$

which may be written as

(5.2) $$\lambda_n = -\sum_{m=1}^{n} \binom{n}{m} \frac{1}{(m-1)!} \sum_{k=1}^{\infty} \frac{(-1)^k}{k} \frac{d^m}{ds^m} \left[2\xi(s) - 1\right]^k \Bigg|_{s=1}$$

For example with $n = 1$ we have

$$\lambda_1 = -\sum_{k=1}^{\infty} \frac{(-1)^k}{k} \frac{d}{ds} \left[2\xi(s)-1\right]^k \bigg|_{s=1}$$

and we see that

$$\frac{d}{ds}\left[2\xi(s)-1\right]^k = 2k\left[2\xi(s)-1\right]^{k-1}\xi'(s)$$

This gives us

$$\sum_{k=1}^{\infty} \frac{(-1)^k}{k} \frac{d}{ds}\left[2\xi(s)-1\right]^k \bigg|_{s=1} = 2\sum_{k=1}^{\infty} \frac{(-1)^k}{k} k\delta_{1,k}\xi'(1)$$

and hence we have

(5.3)     $\lambda_1 = 2\xi'(1)$

Then, using Coffey's result [11] that $\xi^{(n)}(1)$ is positive, we see that

(5.4)     $\lambda_1 > 0$

Using (1.7) we see that

$$\frac{d}{ds}\log 2\xi(s) = \frac{\xi'(s)}{\xi(s)} = -\sum_{k=1}^{\infty}(-1)^k \sigma_k (s-1)^{k-1}$$

which gives us

(5.5)     $2\xi'(1) = \sigma_1$

With $n = 2$ we have

$$\lambda_2 = -\sum_{m=1}^{2}\binom{2}{m}\frac{1}{(m-1)!}\sum_{k=1}^{\infty}\frac{(-1)^k}{k}\frac{d^m}{ds^m}\left[2\xi(s)-1\right]^k\bigg|_{s=1}$$

$$= -2\sum_{k=1}^{\infty}\frac{(-1)^k}{k}\frac{d}{ds}\left[2\xi(s)-1\right]^k\bigg|_{s=1} - \sum_{k=1}^{\infty}\frac{(-1)^k}{k}\frac{d^2}{ds^2}\left[2\xi(s)-1\right]^k\bigg|_{s=1}$$

$$= 2\lambda_1 - \sum_{k=1}^{\infty}\frac{(-1)^k}{k}\frac{d^2}{ds^2}\left[2\xi(s)-1\right]^k\bigg|_{s=1}$$

We have

$$\frac{d^2}{ds^2}\left[2\xi(s)-1\right]^k = 4k(k-1)\left[2\xi(s)-1\right]^{k-2}\left[\xi'(s)\right]^2 + 2k\left[2\xi(s)-1\right]^{k-1}\xi''(s)$$

and thus

$$\frac{d^2}{ds^2}\left[2\xi(s)-1\right]^k\bigg|_{s=1} = 4k(k-1)\delta_{2,k}\left[\xi'(1)\right]^2 + 2k\delta_{1,k}\xi''(1)$$

so that

$$\lambda_2 = 2\lambda_1 - \sum_{k=1}^{\infty}\frac{(-1)^k}{k}\left[4k(k-1)\delta_{2,k}\left[\xi'(1)\right]^2 + 2k\delta_{1,k}\xi''(1)\right]$$

$$= 4\xi'(1) - 4\left[\xi'(1)\right]^2 + 2\xi''(1)$$

Hence we have

(5.6) $\qquad \lambda_2 = 4\xi'(1)\left[1-\xi'(1)\right] + 2\xi''(1)$

We know that $1 > \xi'(1)$ and therefore we deduce that

(5.7) $\qquad \lambda_2 > 0$

We now investigate whether the above process can be generalised.

Often in mathematics we look for divine inspiration but we do not usually expect to obtain it from a canonised saint. This is indeed the source of the next remark. The calculation of $\lambda_n$ effectively involves the derivative of a composite function $g(f(t))$ and the general formula for this was discovered by Francesco Faà di Bruno (1825-1888) who was declared a Saint by Pope John Paul II in St. Peter's Square in Rome on the centenary of his death in 1988 [34].

In [24] di Bruno showed that

(5.8) $\qquad \dfrac{d^m}{dt^m}g(f(t)) = \sum \dfrac{m!}{b_1!b_2!...b_m!}g^{(k)}(f(t))\left(\dfrac{f^{(1)}(t)}{1!}\right)^{b_1}\left(\dfrac{f^{(2)}(t)}{2!}\right)^{b_2}...\left(\dfrac{f^{(m)}(t)}{m!}\right)^{b_m}$

where the sum is over all different solutions in non-negative integers $b_1,...,b_m$ of $b_1 + 2b_2 + ... + mb_m = m$, and $k = b_1 + ... + b_m$.

As reported in [24] this may also be expressed in terms of the (exponential) partial Bell polynomials $B_{m,j}$

$$(5.9) \qquad \frac{d^m}{dt^m} g(f(t)) = \sum_{j=0}^{m} g^{(j)}(f(t)) B_{m,j} \left( f^{(1)}(t), f^{(2)}(t), ..., f^{(m-j+1)}(t) \right)$$

where the (exponential) partial Bell polynomials $B_{m,j} = B_{m,j} \left( x_1, x_2, ..., x_{m-j+1} \right)$ are defined by

$$B_{m,j} = B_{m,j} \left( x_1, x_2, ..., x_{m-j+1} \right) = \sum \frac{m!}{b_1! b_2! ... b_m!} \left( \frac{x_1}{1!} \right)^{b_1} \left( \frac{x_2}{2!} \right)^{b_2} ... \left( \frac{x_1}{m!} \right)^{b_m}$$

where the sum is over all different solutions in non-negative integers $b_1, ..., b_m$ of $b_1 + 2b_2 + ... + mb_m = m$, and $j = b_1 + ... + b_m$.

We now consider $\frac{d^m}{ds^m} \left[ 2\xi(s) - 1 \right]^k$ so that

$$g(x) = x^k$$

$$f(s) = 2\xi(s) - 1$$

$$g^{(j)}(x) = k(k-1)...(k-j+1)x^{k-j}$$

$$g^{(j)}(0) = \delta_{j,k} k!$$

$$f^{(j)}(s) = 2\xi^{(j)}(s)$$

Using (5.9) gives us

$$\frac{d^m}{ds^m} \left[ 2\xi(s) - 1 \right]^k \Big|_{s=1} = \sum_{j=0}^{m} \delta_{j,k} k! B_{m,j} \left( 2\xi^{(1)}(1), 2\xi^{(2)}(1), ..., 2\xi^{(m-j+1)}(1) \right)$$

$$= k! B_{m,k} \left( 2\xi^{(1)}(1), 2\xi^{(2)}(1), ..., 2\xi^{(m-k+1)}(1) \right)$$

Therefore we obtain from (5.2)

$$\lambda_n = -\sum_{m=1}^{n} \binom{n}{m} \frac{1}{(m-1)!} \sum_{k=1}^{\infty} \frac{(-1)^k}{k} k! B_{m,k} \left( 2\xi^{(1)}(1), 2\xi^{(2)}(1), ..., 2\xi^{(m-k+1)}(1) \right)$$

and, since $B_{m,k} = 0$ for $k \geq m+1$, this becomes the finite double summation

$$(5.10) \qquad \lambda_n = \sum_{m=1}^{n} \binom{n}{m} \frac{1}{(m-1)!} \sum_{k=1}^{m} (-1)^{k-1} (k-1)! B_{m,k} \left( 2\xi^{(1)}(1), 2\xi^{(2)}(1), ..., 2\xi^{(m-k+1)}(1) \right)$$

We will see in (7.4) that

$$2\xi(s) = \sum_{n=0}^{\infty} Y_n\left(-\sigma_1, -1!\sigma_2, \ldots, -(n-1)!\sigma_n\right)\frac{s^n}{n!}$$

Hence it follows that

$$2(-1)^m \xi^{(m)}(1) = 2\xi^{(m)}(0) = Y_m\left(-\sigma_1, -1!\sigma_2, \ldots, -(m-1)!\sigma_m\right)$$

We refer to the following inversion relation of Chou et al. [9]

$$(5.11) \quad Y_n = y_n = \sum_{k=1}^{n} B_{n,k}\left(x_1, x_2, \ldots, x_{n-k+1}\right) \Leftrightarrow x_n = \sum_{k=1}^{n}(-1)^{k-1}(k-1)! B_{n,k}\left(y_1, y_2, \ldots, y_{n-k+1}\right)$$

We have seen that

$$2(-1)^n \xi^{(n)}(1) = Y_n\left(-\sigma_1, -1!\sigma_2, \ldots, -(n-1)!\sigma_n\right) = \sum_{k=1}^{n} B_{n,k}\left(x_1, x_2, \ldots, x_{n-k+1}\right)$$

where $x_j = -(j-1)!\sigma_j$. Hence we deduce that

$$-(n-1)!\sigma_n = \sum_{k=1}^{n}(-1)^{k-1}(k-1)! B_{n,k}\left(-2\xi^{(1)}(1), 2\xi^{(2)}(1), \ldots, (-1)^{m+1} 2\xi^{(m+1)}(1)\right)$$

We have [8, p.412]

$$B_{n,k}\left(-2\xi^{(1)}(1), 2\xi^{(2)}(1), \ldots, (-1)^{m+1} 2\xi^{(n+1)}(1)\right) = (-1)^n B_{n,k}\left(2\xi^{(1)}(1), 2\xi^{(2)}(1), \ldots, 2\xi^{(n-k+1)}(1)\right)$$

and using (5.10) above

$$\lambda_n = \sum_{m=1}^{n}\binom{n}{m}\frac{1}{(m-1)!}\sum_{k=1}^{m}(-1)^{k-1}(k-1)! B_{m,k}\left(2\xi^{(1)}(1), 2\xi^{(2)}(1), \ldots, 2\xi^{(m-k+1)}(1)\right)$$

we simply obtain (1.8) again, i.e.

$$\lambda_n = -\sum_{m=1}^{n}\binom{n}{m}(-1)^m \sigma_m$$

$\square$

The logarithmic partition polynomials $L_m$ satisfy the following relations [8, p.424]

$$(5.12) \qquad L_m\left(g_1, g_2, ..., g_m\right) = \sum_{k=1}^{m} (-1)^{k-1} (k-1)! \, B_{m,k}\left(g_1, g_2, ...\right) \quad \text{for } m \geq 1$$

$$L_0 = 0$$

$$(5.13) \qquad \log\left(\sum_{n=0}^{\infty} g_n \frac{s^n}{n!}\right) = \sum_{n=1}^{\infty} L_n \frac{s^n}{n!} \qquad g_0 = 1$$

and referring to (5.10) we see that

$$\lambda_n = \sum_{m=1}^{n} \binom{n}{m} \frac{1}{(m-1)!} L_m\left(2\xi^{(1)}(1), 2\xi^{(2)}(1), ..., 2\xi^{(m)}(1)\right)$$

From (5.13) we have

$$\frac{d^m}{ds^m} \log\left(\sum_{n=0}^{\infty} g_n \frac{s^n}{n!}\right) = \sum_{n=1}^{\infty} L_n n(n-1)...(n-m+1) \frac{s^{n-m}}{n!}$$

and hence

$$\frac{d^m}{ds^m} \log\left(\sum_{n=0}^{\infty} g_n \frac{s^n}{n!}\right)\Bigg|_{s=0} = L_m$$

This then gives us

$$\lambda_n = \sum_{m=1}^{n} \binom{n}{m} \frac{1}{(m-1)!} \frac{d^m}{ds^m} \log\left(\sum_{n=0}^{\infty} g_n \frac{s^n}{n!}\right)\Bigg|_{s=0}$$

where $g_n = 2\xi^{(n)}(1)$.

We see that

$$\sum_{n=0}^{\infty} g_n \frac{s^n}{n!} = 2\sum_{n=0}^{\infty} \xi^{(n)}(1) \frac{s^n}{n!} = 2\sum_{n=0}^{\infty} \xi^{(n)}(0) \frac{(-s)^n}{n!} = 2\xi(-s)$$

and hence we have

$$\lambda_n = \sum_{m=1}^{n} \binom{n}{m} \frac{1}{(m-1)!} \frac{d^m}{ds^m} \log 2\xi(-s)\Bigg|_{s=0}$$

Referring to (2.8) we see that this is equivalent to (1.4.1)

$$\lambda_n = \sum_{m=1}^{n} \binom{n}{m} \frac{1}{(m-1)!} \cdot \frac{d^m}{dp^m} \log 2\xi(s) \bigg|_{s=1}$$

and hence we have simply come full circle!

$\square$

Equation (5.10) may be reindexed to

$$\lambda_n = \sum_{m=1}^{n} \binom{n}{m} \frac{1}{(m-1)!} \sum_{k=0}^{m-1} (-1)^k \, k! \, B_{m,k+1}\left(2\xi^{(1)}(1), 2\xi^{(2)}(1), ..., 2\xi^{(m-k)}(1)\right)$$

and, for example, we have

$$\lambda_1 = B_{1,1}\left(2\xi^{(1)}(1)\right)$$

Since

$$B_{1,1} = x_1$$

we obtain as before

$$\lambda_1 = 2\xi^{(1)}(1)$$

With $n = 2$ we have

$$\lambda_2 = \sum_{m=1}^{2} \binom{2}{m} \frac{1}{(m-1)!} \sum_{k=0}^{m-1} (-1)^k \, k! \, B_{m,k+1}\left(2\xi^{(1)}(1), 2\xi^{(2)}(1), ..., 2\xi^{(m-k)}(1)\right)$$

$$\lambda_2 = 2B_{1,1} + \sum_{k=0}^{1} (-1)^k \, k! \, B_{2,k+1}\left(2\xi^{(1)}(1), 2\xi^{(2)}(1), ..., 2\xi^{(m-k)}(1)\right)$$

$$\lambda_2 = 2B_{1,1} + B_{2,1} - B_{2,2} \qquad\qquad x_n = 2\xi^{(n)}(1)$$

$$B_{2,1} = x_2 \qquad\qquad B_{2,2} = x_1^2$$

$$\lambda_2 = 2x_1 + x_2 - x_1^2$$

and thus we get

$$\lambda_2 = 4\xi^{(1)}(1) + 2\xi''(1) - 4\left[\xi'(1)\right]^2$$

We have

$$\lambda_3 = \sum_{m=1}^{3} \binom{3}{m} \frac{1}{(m-1)!} \sum_{k=0}^{m-1} (-1)^k \, k! \, B_{m,k+1}\left(2\xi^{(1)}(1), 2\xi^{(2)}(1), ..., 2\xi^{(m-k)}(1)\right)$$

$$= \sum_{k=0}^{0} (-1)^k \, k! \, B_{1,k+1} + 3\sum_{k=0}^{1} (-1)^k \, k! \, B_{2,k+1} + 3\sum_{k=0}^{2} (-1)^k \, k! \, B_{3,k+1}$$

$$\lambda_3 = B_{1,1} + 3B_{2,1} - 3B_{2,2} + 3B_{3,1} - 3B_{3,2} + 6B_{3,3}$$

And using

$$B_{3,1} = x_3 \qquad B_{3,2} = 3x_1 x_2 \qquad B_{3,3} = x_1^3$$

we obtain

$$\lambda_3 = x_1 + 3x_2 - 3x_1^2 + 3x_3 - 9x_1 x_2 + 6x_1^3$$

or

$$\lambda_3 = 2\xi^{(1)}(1) + 6\xi^{(2)}(1) - 12\left[\xi^{(1)}(1)\right]^2 + 6\xi^{(3)}(1) - 18\xi^{(1)}(1)\xi^{(2)}(1) + 48\left[\xi^{(1)}(1)\right]^3$$

$\square$

Since

$$\frac{d^m}{ds^m} \log[(s-1)\varsigma(s)] = \frac{d^m}{ds^m} \log[2(s-1)\varsigma(s)]$$

we may also include a factor of $\log 2$ so that

$$\lambda_n = \frac{n}{2}(\gamma + \log \pi) + \sum_{m=2}^{n} \binom{n}{m} 2^{-m} \varsigma(m) + \sum_{m=1}^{n} (-1)^m \binom{n}{m} \frac{1}{(m-1)!} \frac{d^m}{ds^m} \log[2(s-1)\varsigma(s)]\Bigg|_{s=0}$$

Since

$$\frac{d^m}{ds^m} \log[1 + 2(s-1)\varsigma(s) - 1] = -\frac{d^m}{ds^m} \sum_{k=1}^{\infty} \frac{(-1)^k}{k} \left[2(s-1)\varsigma(s) - 1\right]^k$$

$$= -\sum_{k=1}^{\infty} \frac{(-1)^k}{k} \frac{d^m}{ds^m} \left[2(s-1)\varsigma(s) - 1\right]^k$$

we obtain

$$\lambda_n = \frac{n}{2}(\gamma + \log \pi) + \sum_{m=2}^{n} \binom{n}{m} 2^{-m} \varsigma(m) - \sum_{m=1}^{n} (-1)^m \binom{n}{m} \frac{1}{(m-1)!} \sum_{k=1}^{\infty} \frac{(-1)^k}{k} \frac{d^m}{ds^m} [2(s-1)\varsigma(s)-1]^k \bigg|_{s=0}$$

Including the factor of $\log 2$ is convenient because

$$\lim_{s \to 0} \big[ 2(s-1)\varsigma(s) - 1 \big] = 0$$

We now consider $\dfrac{d^m}{ds^m} \big[ 2(s-1)\varsigma(s)-1 \big]^k$ and designating the composite function

$$\big[ 2(s-1)\varsigma(s)-1 \big]^k \equiv g(f(s))$$

so that

$$g(x) = x^k$$

$$f(s) = 2(s-1)\varsigma(s) - 1$$

Differentiation gives us

$$g^{(j)}(x) = k(k-1)\ldots(k-j+1)x^{k-j}$$

$$g^{(j)}(0) = \delta_{j,k} k!$$

$$f^{(j)}(s) = 2 \frac{d^j}{ds^j}(s-1)\varsigma(s)$$

$$\frac{d^j}{ds^j}(s-1)\varsigma(s) = (s-1)\varsigma^{(j)}(s) + j\varsigma^{(j-1)}(s)$$

We then have using (5.9)

$$\frac{d^m}{ds^m} \big[ 2(s-1)\varsigma(s)-1 \big]^k \bigg|_{s=0} = \sum_{j=0}^{m} \delta_{j,k} k! B_{m,j}\left( f^{(1)}(0), f^{(2)}(0), \ldots, f^{(m-j+1)}(0) \right)$$

$$= k! B_{m,k}\left( f^{(1)}(0), f^{(2)}(0), \ldots, f^{(m-k+1)}(0) \right)$$

and from the logarithmic expansion we see that

$$\frac{d^m}{ds^m}\log[(s-1)\varsigma(s)]\Bigg|_{s=0}=-\sum_{k=1}^{\infty}\frac{(-1)^k}{k}k!B_{m,k}\left(f^{(1)}(0),f^{(2)}(0),...,f^{(m-k+1)}(0)\right)$$

Since $B_{m,k}=0$ for $k\geq m+1$ this becomes the finite double summation

$$\frac{d^m}{ds^m}\log[(s-1)\varsigma(s)]\Bigg|_{s=0}=\sum_{k=1}^{m}(-1)^{k-1}(k-1)!B_{m,k}\left(f^{(1)}(0),f^{(2)}(0),...,f^{(m-k+1)}(0)\right)$$

We then obtain

(5.14)    $\lambda_n=\dfrac{n}{2}(\gamma+\log\pi)+\displaystyle\sum_{m=2}^{n}\binom{n}{m}2^{-m}\varsigma(m)$

$$+\sum_{m=1}^{n}(-1)^m\binom{n}{m}\frac{1}{(m-1)!}\sum_{k=1}^{m}(-1)^{k-1}(k-1)!B_{m,k}\left(f^{(1)}(0),f^{(2)}(0),...,f^{(m-k+1)}(0)\right)$$

where

$$f^{(j)}(s)=2\left[j\varsigma^{(j-1)}(0)-\varsigma^{(j)}(0)\right]$$

For example we have

$$\lambda_1=\frac{1}{2}(\gamma+\log\pi)-B_{1,1}\left(f^{(1)}(0)\right)$$

and since $B_{1,1}=x_1$ we obtain again

$$\lambda_1=\frac{1}{2}(\gamma+\log\pi)-2\left[\varsigma(0)-\varsigma^{(1)}(0)\right]$$

We recall (2.11)

$$2\left[m\varsigma^{(m-1)}(0)-\varsigma^{(m)}(0)\right]=Y_m\left(0!b_0,1!b_1,...,(m-1)!b_{m-1}\right)$$

which gives us

$$\sum_{k=1}^{m}(-1)^{k-1}(k-1)!B_{m,k}\left(f^{(1)}(0),f^{(2)}(0),...,f^{(m-k+1)}(0)\right)=\sum_{k=1}^{m}(-1)^{k-1}(k-1)!B_{m,k}\left(Y_1,...,Y_{m-k+1}\right)$$

Noting the inversion relation of Chou et al. [9]

$$y_m = \sum_{k=1}^m B_{m,k}\left(x_1, x_2, ..., x_{m-k+1}\right) \Leftrightarrow x_m = \sum_{k=1}^m (-1)^{k-1}(k-1)! B_{m,k}\left(y_1, y_2, ..., y_{m-k+1}\right)$$

we see that

$$x_m = (m-1)! b_{m-1}$$

and hence we simply recover (2.4).

## 6. Relative magnitudes of $\xi^{(n)}(1)$

Riemann [19] showed that

$$(6.1) \qquad 2\xi(s) = 1 + s(s-1)\int_1^\infty \sum_{k=1}^\infty e^{-\pi k^2 x}\, \frac{x^{(1-s)/2} + x^{s/2}}{x}\, dx$$

and differentiation results in

$$2\xi'(s) = s(s-1)\int_1^\infty \sum_{k=1}^\infty e^{-\pi k^2 x}\, \frac{\log \sqrt{x}}{x}\left[(-1)x^{(1-s)/2} + x^{s/2}\right] dx + (2s-1)\int_1^\infty \sum_{k=1}^\infty e^{-\pi k^2 x}\, \frac{x^{(1-s)/2} + x^{s/2}}{x}\, dx$$

With $s = 1$ we have

$$2\xi'(1) = \int_1^\infty \sum_{k=1}^\infty e^{-\pi k^2 x}\, \frac{1 + x^{1/2}}{x}\, dx$$

and we see that $\xi'(1) > 0$. A further differentiation gives us

$$2\xi''(s) = s(s-1)\int_1^\infty \sum_{k=1}^\infty e^{-\pi k^2 x}\, \frac{\log^2 \sqrt{x}}{x}\left[x^{(1-s)/2} + x^{s/2}\right] dx$$

$$+ 2(2s-1)\int_1^\infty \sum_{k=1}^\infty e^{-\pi k^2 x}\, \frac{\log \sqrt{x}}{x}\left[(-1)x^{(1-s)/2} + x^{s/2}\right] dx + 2\int_1^\infty \sum_{k=1}^\infty e^{-\pi k^2 x}\, \frac{x^{(1-s)/2} + x^{s/2}}{x}\, dx$$

and we have

$$2\xi''(1) = 2\int_1^\infty \sum_{k=1}^\infty e^{-\pi k^2 x}\, \frac{\log \sqrt{x}}{x}\left[x^{1/2} - 1\right] dx + 2\int_1^\infty \sum_{k=1}^\infty e^{-\pi k^2 x}\, \frac{1 + x^{1/2}}{x}\, dx$$

$$= 2\int\limits_1^\infty \sum_{k=1}^\infty e^{-\pi k^2 x} \frac{\log \sqrt{x}}{x} \left[ x^{1/2} - 1 \right] dx + 4\xi'(1)$$

Therefore we see that

(6.2)    $\xi''(1) > 2\,\xi'(1)$

We write (6.1) as

$$2\xi(s) = 1 + s(s-1) f(s)$$

where

$$f(s) = \int\limits_1^\infty \sum_{k=1}^\infty e^{-\pi k^2 x} \frac{x^{(1-s)/2} + x^{s/2}}{x} dx$$

Using the Leibniz rule we have

$$2\xi^{(n)}(s) = \sum_{j=0}^n \binom{n}{j} f^{(n-j)}(s) \frac{d^{(j)}}{ds^{(j)}} s(s-1)$$

which immediately simplifies for $n \geq 2$ to

$$2\xi^{(n)}(s) = s(s-1)\binom{n}{0} f^{(n)}(s) + (2s-1)\binom{n}{1} f^{(n-1)}(s) + 2\binom{n}{2} f^{(n-2)}(s)$$

With $s = 0$ we obtain

$$\xi^{(n)}(0) = \binom{n}{2} f^{(n-2)}(0) - \frac{1}{2}\binom{n}{1} f^{(n-1)}(0)$$

We see that

$$f^{(j)}(s) = \frac{1}{2^j} \int\limits_1^\infty \sum_{k=1}^\infty e^{-\pi k^2 x} \frac{\log^j x}{x} \left[ (-1)^j x^{(1-s)/2} + x^{s/2} \right] dx$$

$$f^{(j)}(0) = \frac{1}{2^j} \int\limits_1^\infty \sum_{k=1}^\infty e^{-\pi k^2 x} \frac{\log^j x}{x} \left[ (-1)^j x^{1/2} + 1 \right] dx$$

$$= \frac{(-1)^j}{2^j} \int\limits_1^\infty \sum_{k=1}^\infty e^{-\pi k^2 x} \frac{\log^j x}{x} \left[ x^{1/2} + (-1)^j \right] dx$$

Therefore we have

$$\xi^{(n)}(0) = (-1)^n \frac{n}{2} \int_1^\infty \sum_{k=1}^\infty e^{-\pi k^2 x} \frac{\log^{n-2}\sqrt{x}}{x} \left( (n-1)\left[\sqrt{x} + (-1)^n\right] + \left[\sqrt{x} - (-1)^n\right]\log\sqrt{x} \right) dx$$

Since $\xi^{(n)}(0) = (-1)^n \xi^{(n)}(1)$ we have

$$(6.3) \qquad \xi^{(n)}(1) = \frac{n}{2} \int_1^\infty \sum_{k=1}^\infty e^{-\pi k^2 x} \frac{\log^{n-2}\sqrt{x}}{x} \left( (n-1)\left[\sqrt{x} + (-1)^n\right] + \left[\sqrt{x} - (-1)^n\right]\log\sqrt{x} \right) dx$$

and, since $\sqrt{x} \pm (-1)^n \geq 0$ for $x \geq 1$, we see that for $n \geq 2$

$$(6.4) \qquad \xi^{(n)}(1) > 0$$

As reported by Keiper [25] we have

$$2\xi(s) = \sum_{n=0}^\infty \alpha_n (s-1)^n = 2\sum_{n=0}^\infty \frac{\xi^{(n)}(1)}{n!}(s-1)^n$$

where

$$\alpha_0 = 1 \qquad\qquad \alpha_n = \beta_{n-2} + \beta_{n-1} \quad \text{for } n \geq 1$$

$$\beta_{-1} = 0 \qquad\qquad \beta_0 = 1 + \frac{1}{2}\gamma - \log\left(2\sqrt{\pi}\right)$$

$$\beta_n = \frac{1}{n!} \int_1^\infty \sum_{k=1}^\infty e^{-\pi k^2 x} \frac{\log^n \sqrt{x}}{x}\left(\sqrt{x} + (-1)^n\right) dx \quad \text{for } n \geq 1$$

$$\xi^{(n)}(1) = \frac{1}{2}\alpha_n n!$$

$$\alpha_0 = 1$$

$$\alpha_1 = 1 + \frac{1}{2}\gamma - \log\left(2\sqrt{\pi}\right)$$

$$\alpha_2 = 1 + \frac{1}{2}\gamma - \log\left(2\sqrt{\pi}\right) + \beta_1$$

$$\alpha_3 = \beta_1 + \beta_2$$

It is easily seen that $\beta_n > 0$ for $n \geq 1$ which implies that $\alpha_n > 0$ for $n \geq 3$. Therefore we see again that $\xi^{(n)}(1) > 0$ for (at least) $n \geq 3$. Curiously, this point was not specifically mentioned in Keiper's paper [25]; this was first proved (and extended) by Coffey [11] in 2004 using a different formulation, namely

$$\xi^{(2m)}(1) = \frac{2\pi}{2^{2m}} \int_1^\infty \sum_{n=1}^\infty \left( \pi n^2 x - \frac{3}{2} \right) n^2 e^{-\pi n^2 x} \sqrt{x} \left[ \sqrt{x} + 1 \right] \log^{2m} x \, dx$$

$$\xi^{(2m+1)}(1) = \frac{\pi}{2^{2m}} \int_1^\infty \sum_{n=1}^\infty \left( \pi n^2 x - \frac{3}{2} \right) n^2 e^{-\pi n^2 x} \sqrt{x} \left[ \sqrt{x} - 1 \right] \log^{2m+1} x \, dx$$

Coffey's result [11] may also be obtained by differentiating Riemann's formula [19, p.17]

$$\xi(s) = 4 \int_1^\infty \frac{d}{dx} \left[ x^{3/2} \omega'(x) \right] x^{-1/4} \cosh\left[ \frac{1}{2}\left( s - \frac{1}{2} \right) \log x \right] dx$$

where $\omega(x) = \sum_{k=1}^\infty e^{-\pi k^2 x}$.

We immediately see that for $n \geq 3$

$$\frac{2}{n!} \left[ \xi^{(n+1)}(1) - \xi^{(n)}(1) \right] = n\beta_{n-1} + (n+1)\beta_n - \beta_{n-2}$$

so that

$$\frac{2}{n!} \left[ \xi^{(n+1)}(1) - \xi^{(n)}(1) \right] =$$

$$\frac{1}{n!} \int_1^\infty \sum_{k=1}^\infty e^{-\pi k^2 x} \frac{\log^{n-2} \sqrt{x}}{x} \left[ n^2 \log \sqrt{x} \left( \sqrt{x} - (-1)^n \right) + (n+1) \log^2 \sqrt{x} \left( \sqrt{x} + (-1)^n \right) - n(n-1) \left( \sqrt{x} + (-1)^n \right) \right] dx$$

With the substitution $t = \sqrt{x}$ the integral becomes

$$= \frac{2}{n!} \int_1^\infty \sum_{k=1}^\infty e^{-\pi k^2 t^2} \frac{\log^{n-2} t}{t} \left[ n^2 \log t \left( t - (-1)^n \right) + (n+1) \log^2 t \left( t + (-1)^n \right) - n(n-1) \left( t + (-1)^n \right) \right] dt$$

For convenience we designate $f_n(t)$ for $n \geq 3$ as

$$f_n(t) = \left[ n^2 \log t \left( t - (-1)^n \right) + (n+1) \log^2 t \left( t + (-1)^n \right) - n(n-1) \left( t + (-1)^n \right) \right] \log^{n-2} t$$

and we first of all consider the specific case where $n = 3$, giving us

$$f_3(t) = \left[ 9 \log t (t+1) + 4 \log^2 t (t-1) - 6(t-1) \right] \log t$$

Differentiation results in

$$f_3'(t) = \left[ 9 \log t (t+1) + 4 \log^2 t (t-1) - 6(t-1) \right] \frac{1}{t}$$

$$+ \left[ 9 \log t + 9 \frac{1}{t} (t+1) + 4 \log^2 t + 8 \frac{1}{t} \log t (t-1) - 6 \right] \log t$$

$$= \left[ 9 \log t (t+1) + 4 \log^2 t (t-1) - 6(t-1) \right] \frac{1}{t}$$

$$+ \left[ 9 \log t + 9 \frac{1}{t} + 4 \log^2 t + 8 \frac{1}{t} \log t (t-1) + 3 \right] \log t$$

where we note that the second term in parentheses is positive for $t \geq 1$. We now consider the first term

$$g(t) = 9 \log t (t+1) + 4 \log^2 t (t-1) - 6(t-1) \qquad g(1) = 0$$

Differentiation results in

$$g'(t) = 9 \log t + 9 \frac{1}{t} (t+1) + 4 \log^2 t + 8 \frac{1}{t} \log t (t-1) - 6$$

$$= 9 \log t + 9 \frac{1}{t} + 4 \log^2 t + 8 \frac{1}{t} \log t (t-1) + 3$$

and we note that $g'(t)$ is positive for $t \geq 1$. We therefore see that $f_3'(t) \geq 0$ for all $t \geq 1$, thereby proving that $f_3(t)$ is monotonic increasing for $t \geq 1$. Since $f_3(1) = 0$, we deduce that $f_3(t) \geq 0$ for all $t \geq 1$.

We now consider the case where $n \geq 3$

$$f_n'(t) = \left[ n^2 \log t \left( t - (-1)^n \right) + (n+1) \log^2 t \left( t + (-1)^n \right) - n(n-1) \left( t + (-1)^n \right) \right] (n-2) \frac{1}{t} \log^{n-3} t$$

$$+ \left[ n^2 \log t + n^2 \frac{1}{t} \left( t - (-1)^n \right) + (n+1) \log^2 t + 2(n+1) \log t \frac{1}{t} \left( t + (-1)^n \right) - n(n-1) \right] \log^{n-2} t$$

and we write this as

$$f_n'(t) = g(t)(n-2)\frac{1}{t}\log^{n-3} t + h(t)\log^{n-2} t$$

First of all, we see that $h(t)$ is certainly positive when $n$ is odd and $t \geq 1$.

We have

$$g'(t) = n^2 \log t + n^2 \frac{1}{t}\left(t - (-1)^n\right) + (n+1)\log^2 t + (n+1)\frac{1}{t}\left(t + (-1)^n\right) - n(n-1)$$

$$= n^2 \log t - (-1)^n n^2 \frac{1}{t} + (n+1)\log^2 t + (n+1)\frac{1}{t}\left(t + (-1)^n\right) + n$$

and we note that $g'(t)$ is certainly positive when $n$ is odd and $n \geq 3$. Therefore we conclude that $f_n(t) \geq 0$ when $n \geq 3$ and $n$ is odd. Hence we have determined that

(6.5)     $$\xi^{(2m+2)}(1) > \xi^{(2m+1)}(1) \text{ for } m \geq 1$$

so that, for example, we have

(6.5.1)     $$\xi^{(4)}(1) > \xi^{(3)}(1)$$

We shall determine below in (6.6) that

$$\xi^{(2)}(1) > \xi^{(3)}(1)$$

$\square$

Using (6.3) we see that

$$\xi^{(2)}(1) - \xi^{(3)}(1) = \int_1^\infty \sum_{k=1}^\infty e^{-\pi k^2 x} \frac{1}{x}\left(\sqrt{x} + 1 - 2\left[\sqrt{x} - 1\right]\log\sqrt{x} - \frac{3}{2}\left[\sqrt{x} + 1\right]\log^2\sqrt{x}\right)dx$$

With the substitution $t = \sqrt{x}$ the integral becomes

$$= 2\int_1^\infty \sum_{k=1}^\infty e^{-\pi k^2 t^2} \frac{1}{t}\left(t + 1 - 2[t-1]\log t - \frac{3}{2}[t+1]\log^2 t\right)dt$$

With the substitution $u = 1/t$ the integral becomes

$$= 2\int_0^1 \sum_{k=1}^\infty e^{-\pi k^2/u^2} u^2\left(1 + u + 2[1-u]\log u - \frac{3}{2}[1+u]\log^2 u\right)du$$

and we now consider that part of the integrand defined by

$$g(u) = u^2 \left( 1 + u + 2[1-u]\log u - \frac{3}{2}[1+u]\log^2 u \right)$$

where we note that $g(0) = 0$. We have the derivative

$$g'(u) = u^2 + 4u - 9u^2 \log u - 3u \log u - \frac{9}{2} u^2 \log^2 u$$

$$= u^2 + 4u - \frac{9}{2} u^2 \log u [2 + \log u] - 3u \log u$$

If $2 + \log u \geq 0$, then $g'(u) > 0$ because $\log u$ is negative in the interval $(0,1)$. Let us now consider the other possibility where $-\log u > 2$; therefore in the interval $\left[ 0, e^{-2} \right]$ we have the inequality

$$g'(u) > u^2 + 4u + 18u^2 + 6u - \frac{9}{2} u^2 \log^2 u$$

$$= uh(u)$$

where $h(u)$ is defined as

$$h(u) = 10 + 19u - \frac{9}{2} u \log^2 u$$

The smallest value of $h(u)$ arises when (i) $\frac{9}{2} u \log^2 u$ is at its maximum value and (ii) the other component, $10 + 19u$, is at its minimum value. We now consider the function $f(u)$ defined by

$$f(u) = \frac{9}{2} u \log^2 u$$

We note that $f(u) \geq 0$ and that $f(0) = f(1) = 0$ and we have

$$f'(u) = \frac{9}{2} \log u [2 + \log u]$$

It is easily seen that $f'(u) \geq 0$ in the interval $\left(0, e^{-2}\right)$ and $f'(u) \leq 0$ in the interval $\left(e^{-2}, 1\right)$. Furthermore, $f'(u) = 0$ when $u = e^{-2}$. Therefore, the maximum value of $f(u)$ is attained when $u = e^{-2}$, i.e. the maximum value of $f(u)$ is $18e^{-2}$.

Accordingly, the smallest value of $h(u)$ is equal to $10 - 18e^{-2}$ (which is a positive number). Hence we have determined that $h(u) > 0$ and thereby deduce that $g'(u) > 0$. We have therefore proved that $g(u)$ is monotonic increasing and, since $g(0) = 0$, we see that $g(u) > 0$.

Hence, since the integrand is non-negative, we deduce that

(6.6) $\qquad \xi^{(2)}(1) > \xi^{(3)}(1)$

This inequality is employed in Section 7 below.

Further work is required to determine whether similar inequalities exist for other combinations of $m$.

## 7. The sigma constants $\sigma_n$

It is well known that if

(7.1) $\qquad \log h(x) = b_0 + \sum_{n=1}^{\infty} \frac{b_n}{n} x^n$

then we have in terms of the (exponential) complete Bell polynomials (see Section 10 of this paper)

(7.2) $\qquad h(x) = e^{b_0} \sum_{n=0}^{\infty} Y_n\left(b_1, 1!b_2, ..., (n-1)!b_n\right) \frac{x^n}{n!}$

and

(7.3) $\qquad h^{\alpha}(x) = e^{\alpha b_0} \sum_{n=0}^{\infty} Y_n\left(\alpha b_1, 1!\alpha b_2, ..., (n-1)!\alpha b_n\right) \frac{x^n}{n!}$

If we are only dealing with real values of $x$ then (7.1) implies that

(7.3.1) $\qquad h(x) > 0$

Referring to (1.7.2)

$$\log 2\xi(s) = -\sum_{k=1}^{\infty} \frac{\sigma_k}{k} s^k$$

we then determine that

(7.4) $$2\xi(s) = \sum_{n=0}^{\infty} Y_n\left(-\sigma_1, -1!\sigma_2, \ldots, -(n-1)!\sigma_n\right) \frac{s^n}{n!}$$

or equivalently

$$2\xi(1-s) = \sum_{n=0}^{\infty} Y_n\left(-\sigma_1, -1!\sigma_2, \ldots, -(n-1)!\sigma_n\right) \frac{(1-s)^n}{n!}$$

Differentiation gives us

(7.5) $$2\xi^{(n)}(0) = 2(-1)^n \xi^{(n)}(1) = Y_n\left(-\sigma_1, -1!\sigma_2, \ldots, -(n-1)!\sigma_n\right)$$

Since $\xi^{(n)}(1)$ is positive, we see that $Y_n\left(-\sigma_1, -1!\sigma_2, \ldots, -(n-1)!\sigma_n\right)$ has the same sign as $(-1)^n$.

With $n = 1$ we see that

(7.6) $$2\xi^{(1)}(1) = -Y_1\left(-\sigma_1\right) = \sigma_1$$

Hence we see that $\sigma_1 > 0$. By calculation we find that $1 > \sigma_1 > 0$ and hence we have

(7.7) $$\sigma_1 > \sigma_1^2$$

Referring to (1.8)

$$\lambda_n = -\sum_{m=1}^{n} (-1)^m \binom{n}{m} \sigma_m$$

we see that

(7.8) $$\lambda_1 = \sigma_1 > 0$$

With $n = 2$ in (7.5) we obtain

$$2\xi^{(2)}(1) = Y_2\left(-\sigma_1, -1!\sigma_2\right)$$

and thus

(7.9)     $2\xi^{(2)}(1) = \sigma_1^2 - \sigma_2$

This then tells us that

$$\sigma_1^2 > \sigma_2$$

We have from (1.8)

$$\lambda_2 = 2\sigma_1 - \sigma_2$$

$$> 2\sigma_1 - \sigma_1^2$$

and since $-\sigma_1^2 > -\sigma_1$ we see that

(7.10)     $\lambda_2 > 2\sigma_1 - \sigma_1 = \sigma_1$

Hence we have

(7.11)     $\lambda_2 > \lambda_1 > 0$

We have seen in (6.2) that $\xi^{(2)}(1) > 2\xi^{(1)}(1)$ and thus $\xi^{(2)}(1) > \xi^{(1)}(1)$ implies that

$$\sigma_1^2 - \sigma_2 > \sigma_1$$

or equivalently

$$\sigma_1(\sigma_1 - 1) - \sigma_2 > 0$$

Since the first term on the left-hand side is negative, we deduce that $\sigma_2$ must be negative.

Furthermore, since $\lambda_2 = 2\sigma_1 - \sigma_2 > 0$ we determine that

$$2\sigma_1 > \sigma_2$$

but this does not really provide any more useful information (because $\sigma_1$ and $\sigma_2$ are of opposite signs).

We now consider the case where $n = 3$. We have from (7.5)

$$2\xi^{(3)}(1) = -Y_3\left(-\sigma_1, -1!\,\sigma_2, -2!\,\sigma_3\right)$$

and since

$$Y_3(x_1, x_2, x_3) = x_1^3 + 3x_1 x_2 + x_3$$

we obtain

(7.12) $\qquad 2\xi^{(3)}(1) = \sigma_1^3 - 3\sigma_1\sigma_2 + 2\sigma_3$

Since $\xi^{(3)}(1)$ is positive we have

(7.13) $\qquad \sigma_1^3 - 3\sigma_1\sigma_2 + 2\sigma_3 > 0$

and (1.8) gives us

$$\lambda_3 = 3\sigma_1 - 3\sigma_2 + \sigma_3 = 3(\sigma_1 - \sigma_2) + \sigma_3$$

Then using (7.13) we see that

$$\lambda_3 > 3\sigma_1 - 3\sigma_2 + \frac{3}{2}\sigma_1\sigma_2 - \frac{1}{2}\sigma_1^3$$

$$= 3\sigma_1 + 3\sigma_2\left(\frac{1}{2}\sigma_1 - 1\right) - \frac{1}{2}\sigma_1^3$$

$$> \frac{5}{2}\sigma_1 + 3\sigma_2\left(\frac{1}{2}\sigma_1 - 1\right)$$

Therefore, since the second term on the right-hand side is also positive, we conclude that

(7.14) $\qquad \lambda_3 > 0$

Referring back to

$$\lambda_3 \quad > \frac{5}{2}\sigma_1 + 3\sigma_2\left(\frac{1}{2}\sigma_1 - 1\right)$$

$$= \frac{5}{2}\sigma_1 - \sigma_2 + \frac{3}{2}\sigma_2\left(\sigma_1 - \frac{4}{3}\right)$$

$$> 2\sigma_1 - \sigma_2 + \frac{3}{2}\sigma_2\left(\sigma_1 - \frac{4}{3}\right)$$

Hence, since $\lambda_2 = 2\sigma_1 - \sigma_2$, we determine that

(7.15) $\qquad \lambda_3 > \lambda_2$

With reference to (1.8) we have the binomial transform

$$\sigma_n = -\sum_{m=1}^{n} (-1)^m \binom{n}{m} \lambda_m$$

and in particular we have

$$\sigma_2 = 2\lambda_1 - \lambda_2$$

Since $\sigma_2$ is negative, this implies that

(7.16) $\qquad \lambda_2 > 2\lambda_1$

Similarly we have

$$\sigma_3 = 3\lambda_1 - 3\lambda_2 + \lambda_3$$

Using (6.6) $\xi^{(2)}(1) > \xi^{(3)}(1)$ we have

$$\sigma_1^2 - \sigma_2 \; > \; \sigma_1^3 - 3\sigma_1\sigma_2 + 2\sigma_3$$

Equation (6.5)

$$\xi^{(2m+2)}(1) \; > \; \xi^{(2m+1)}(1) \; \text{ for } m \geq 1$$

may also be useful in investigating the properties of the higher orders of $\lambda_n$, albeit the algebra will become increasingly tedious.

$\qquad\qquad\qquad\qquad\qquad\qquad\qquad\qquad\qquad\qquad\qquad\qquad\qquad\qquad$ $\square$

Reference to (7.3) and (7.4) shows that

$$\left[2\xi(s)\right]^j = \sum_{n=0}^{\infty} Y_n\left(-j\sigma_1, -1!\, j\sigma_2, \ldots, -(n-1)!\, j\sigma_n\right)\frac{s^n}{n!}$$

and using (1.3) $\xi(s) = \xi(1-s)$ this becomes

$$= \sum_{n=0}^{\infty} Y_n \left( -j\sigma_1, -1! \, j\sigma_2, ..., -(n-1)! \, j\sigma_n \right) \frac{(1-s)^n}{n!}$$

The binomial theorem gives us

$$[2\xi(s)-1]^k = \sum_{j=0}^{k} \binom{k}{j}(-1)^{k-j}[2\xi(s)]^j$$

$$= \sum_{j=0}^{k} \binom{k}{j}(-1)^{k-j} \sum_{n=0}^{\infty} Y_n \left( -j\sigma_1, -1! \, j\sigma_2, ..., -(n-1)! \, j\sigma_n \right) \frac{(1-s)^n}{n!}$$

We see that

$$\frac{d^m}{ds^m}[2\xi(s)-1]^k$$

$$= (-1)^m \sum_{j=0}^{k} \binom{k}{j}(-1)^{k-j} \sum_{n=0}^{\infty} Y_n \left( -j\sigma_1, -1! \, j\sigma_2, ..., -(n-1)! \, j\sigma_n \right) n(n-1)...(n-m+1) \frac{(1-s)^{n-m}}{n!}$$

and thus

$$\frac{d^m}{ds^m}[2\xi(s)-1]^k \bigg|_{s=1} = (-1)^m \sum_{j=0}^{k} \binom{k}{j}(-1)^{k-j} Y_m \left( -j\sigma_1, -1! \, j\sigma_2, ..., -(m-1)! \, j\sigma_m \right)$$

We recall (5.2)

$$\lambda_n = -\sum_{m=1}^{n} \binom{n}{m} \frac{1}{(m-1)!} \sum_{k=1}^{\infty} \frac{(-1)^k}{k} \frac{d^m}{ds^m}[2\xi(s)-1]^k \bigg|_{s=1}$$

which gives us

$$(7.17) \quad \lambda_n = \sum_{m=1}^{n} \binom{n}{m} \frac{(-1)^{m+1}}{(m-1)!} \sum_{k=1}^{\infty} \frac{1}{k} \sum_{j=0}^{k} \binom{k}{j}(-1)^j Y_m \left( -j\sigma_1, -1! \, j\sigma_2, ..., -(m-1)! \, j\sigma_m \right)$$

For $n=1$ we have

$$\lambda_1 = \sum_{k=1}^{\infty} \frac{1}{k} \sum_{j=0}^{k} \binom{k}{j}(-1)^j Y_1 \left( -j\sigma_1 \right)$$

$$= -\sigma_1 \sum_{k=1}^{\infty} \frac{1}{k} \sum_{j=0}^{k} \binom{k}{j}(-1)^j j$$

$$= \sigma_1 \sum_{k=1}^{\infty} \delta_{1,k}$$

where we used the identity

$$\sum_{j=0}^{k} \binom{k}{j} (-1)^j j = -k \delta_{1,k}$$

This is easily derived using the binomial theorem

$$(1-x)^k = \sum_{j=0}^{k} \binom{k}{j} (-1)^j x^j$$

whereby differentiation results in

$$(7.18) \qquad -k(1-x)^{k-1} = \sum_{j=0}^{k} \binom{k}{j} (-1)^j j x^{j-1}$$

and the result follows by letting $x = 1$.

Hence we obtain

$$\lambda_1 = \sigma_1$$

For $n = 2$ we have

$$\lambda_2 = \sum_{m=1}^{2} \binom{2}{m} \frac{(-1)^{m+1}}{(m-1)!} \sum_{k=1}^{\infty} \frac{1}{k} \sum_{j=0}^{k} \binom{k}{j} (-1)^j Y_m \left( -j\sigma_1, -1! \, j\sigma_2, \ldots, -(m-1)! \, j\sigma_m \right)$$

$$= 2 \sum_{k=1}^{\infty} \frac{1}{k} \sum_{j=0}^{k} \binom{k}{j} (-1)^j Y_1 \left( -j\sigma_1 \right) - \sum_{k=1}^{\infty} \frac{1}{k} \sum_{j=0}^{k} \binom{k}{j} (-1)^j Y_2 \left( -j\sigma_1, -1! \, j\sigma_2 \right)$$

$$= 2\sigma_1 - \sum_{k=1}^{\infty} \frac{1}{k} \sum_{j=0}^{k} \binom{k}{j} (-1)^j \left[ j^2 \sigma_1^2 - j\sigma_2 \right]$$

$$= 2\sigma_1 - \sigma_1^2 \sum_{k=1}^{\infty} \frac{1}{k} \sum_{j=0}^{k} \binom{k}{j} (-1)^j j^2 + \sigma_2 \sum_{k=1}^{\infty} \frac{1}{k} \sum_{j=0}^{k} \binom{k}{j} (-1)^j j$$

$$= 2\sigma_1 - \sigma_1^2 \sum_{k=1}^{\infty} \frac{1}{k} \sum_{j=0}^{k} \binom{k}{j} (-1)^j j^2 - \sigma_2$$

We multiply (7.18) by $x$ and differentiate to obtain

$$k(k-1)x(1-x)^{k-2} - k(1-x)^{k-1} = \sum_{j=0}^{k} \binom{k}{j}(-1)^j j^2 x^{j-1}$$

and we see that

$$\sum_{j=0}^{k} \binom{k}{j}(-1)^j j^2 = k(k-1)\delta_{2,k} - k\delta_{1,k}$$

Therefore we have

$$\sum_{k=1}^{\infty} \frac{1}{k}\sum_{j=0}^{k} \binom{k}{j}(-1)^j j^2 = \sum_{k=1}^{\infty} \frac{1}{k}[k(k-1)\delta_{2,k} - k\delta_{1,k}] = 0$$

and hence we obtain

$$\lambda_2 = 2\sigma_1 - \sigma_2$$

It is clear that the calculations may be extended to higher orders of $\lambda_n$; however I suspect that this will simply lead to (1.8)

$$\lambda_n = -\sum_{j=1}^{n} (-1)^j \binom{n}{j}\sigma_j$$

We may also note from [21] that the Bernoulli polynomials are given by

$$B_p(u) = \sum_{n=0}^{\infty} \frac{1}{n+1}\sum_{k=0}^{n} \binom{n}{k}(-1)^k (u+k)^p$$

and hence the Bernoulli numbers are given by

$$B_p = \sum_{n=0}^{\infty} \frac{1}{n+1}\sum_{k=0}^{n} \binom{n}{k}(-1)^k k^p$$

I initially thought that these numbers would feature in the above analysis but, as we shall see below, this is not the case.

$\square$

We now refer to (10.18)

$$Y_m(\alpha x_1, ..., \alpha x_m) = \sum_{l=1}^{m} \alpha^l B_{m,l}(x_1, ..., x_{m-l+1})$$

which shows that we may write (7.17) as

$$\lambda_n = \sum_{m=1}^n \binom{n}{m} \frac{(-1)^{m+1}}{(m-1)!} \sum_{k=1}^\infty \frac{1}{k} \sum_{j=0}^k \binom{k}{j} (-1)^j \sum_{l=1}^m j^l B_{m,l} \left(-\sigma_1, -1!\sigma_2, ..., -(m-l)!\sigma_{m-l+1}\right)$$

We see that

$$\sum_{j=0}^k \binom{k}{j}(-1)^j \sum_{l=1}^m j^l B_{m,l} = \sum_{l=1}^m B_{m,l} \sum_{j=0}^k \binom{k}{j}(-1)^j j^l$$

and hence we have

$$\lambda_n = \sum_{m=1}^n \binom{n}{m} \frac{(-1)^{m+1}}{(m-1)!} \sum_{l=1}^m B_{m,l} \sum_{k=1}^\infty \frac{1}{k} \sum_{j=0}^k \binom{k}{j}(-1)^j j^l$$

We have the well known representation for the Stirling numbers of the second kind [8, p.289]

$$(-1)^l S(l,k)k! = \sum_{j=0}^k \binom{k}{j}(-1)^j j^l$$

and thus

$$\lambda_n = \sum_{m=1}^n \binom{n}{m} \frac{(-1)^{m+1}}{(m-1)!} \sum_{l=1}^m B_{m,l}(-1)^l \sum_{k=1}^\infty \frac{1}{k} S(l,k)k!$$

Since $S(l,k) = 0$ for $k \geq l+1$ we have

$$(7.19) \quad \lambda_n = \sum_{m=1}^n \binom{n}{m} \frac{(-1)^{m+1}}{(m-1)!} \sum_{l=1}^m B_{m,l}(-1)^l \sum_{k=1}^l (k-1)! S(l,k)$$

and, for example, this gives us

$$\lambda_1 = -B_{1,1}S(1,1) = -B_{1,1} = \sigma_1$$

$\square$

We note the recurrence relation [8, p.448]

$$Y_n(x_1 + y_1, ..., x_n + y_n) = \sum_{k=0}^n \binom{n}{k} Y_{n-k}(x_1, ..., x_{n-k}) Y_k(y_1, ..., y_k)$$

which shows us that

$$Y_n(-2x_1, \ldots, -2x_n) = \sum_{k=0}^{n} \binom{n}{k} Y_{n-k}(-x_1, \ldots, -x_{n-k}) Y_k(-x_1, \ldots, -y_k)$$

and this may be generalised to the Bell polynomials with the following arguments

$$Y_m\left(-j\sigma_1, -1!\, j\sigma_2, \ldots, -(m-1)!\, j\sigma_m\right)$$

## 8. Further applications of the Bell polynomials

We recall (1.12)

$$\log[(s-1)\varsigma(s)] = -\sum_{k=1}^{\infty} \frac{\eta_{k-1}}{k}(s-1)^k$$

and with $s-1 \rightarrow s$ this becomes

$$\log[s\varsigma(s+1)] = -\sum_{k=1}^{\infty} \frac{\eta_{k-1}}{k} s^k$$

Applying (7.2) gives us

$$s\varsigma(s+1) = \sum_{n=0}^{\infty} Y_n\left(-\eta_0, -1!\,\eta_1, \ldots, -(n-1)!\,\eta_{n-1}\right) \frac{s^n}{n!}$$

and with $s-1 \rightarrow s$ we obtain

$$(s-1)\varsigma(s) = \sum_{n=0}^{\infty} Y_n\left(-\eta_0, -1!\,\eta_1, \ldots, -(n-1)!\,\eta_{n-1}\right) \frac{(s-1)^n}{n!}$$

so that

$$\frac{d^m}{ds^m}[(s-1)\varsigma(s)]\bigg|_{s=1} = Y_m\left(-\eta_0, -1!\,\eta_1, \ldots, -(m-1)!\,\eta_{m-1}\right)$$

Using (1.19) we see that

$$\frac{d^m}{ds^m}[(s-1)\varsigma(s)]\bigg|_{s=1} = (-1)^{m-1} m\gamma_{m-1}$$

and hence we have

(8.1)     $(-1)^{m-1} m\gamma_{m-1} = Y_m\left(-\eta_0, -1!\,\eta_1, \ldots, -(m-1)!\,\eta_{m-1}\right)$

which we originally reported in [17].

Referring to (3.2), i.e. $\eta_k = (-1)^{k+1}\varepsilon_k$ we obtain

$$(-1)^{m-1} m \gamma_{m-1} = Y_m \left( \varepsilon_0, -1!\varepsilon_1, ..., (-1)^{m+1}(m-1)!\varepsilon_{m-1} \right)$$

and using (10.6) this becomes

$$= (-1)^m Y_m \left( -\varepsilon_0, -1!\varepsilon_1, ..., -(m-1)!\varepsilon_{m-1} \right)$$

However, no discernible sign pattern for $\gamma_{m-1}$ emerges from this representation.

We also see that

$$\log \frac{1}{(s-1)\varsigma(s)} = \sum_{k=1}^{\infty} \frac{\eta_{k-1}}{k} (s-1)^k$$

so that

$$\log \frac{1}{s\varsigma(s+1)} = \sum_{k=1}^{\infty} \frac{\eta_{k-1}}{k} s^k$$

Therefore applying (7.2) gives us

$$\frac{1}{s\varsigma(s+1)} = \sum_{n=0}^{\infty} Y_n \left( \eta_0, 1!\eta_1, ..., (n-1)!\eta_{n-1} \right) \frac{s^n}{n!}$$

so that

$$\frac{1}{(s-1)\varsigma(s)} = \sum_{n=0}^{\infty} Y_n \left( \eta_0, 1!\eta_1, ..., (n-1)!\eta_{n-1} \right) \frac{(s-1)^n}{n!}$$

We then obtain

$$\frac{d^m}{ds^m} \frac{1}{(s-1)\varsigma(s)} \bigg|_{s=1} = Y_m \left( \eta_0, 1!\eta_1, ..., (m-1)!\eta_{m-1} \right)$$

Employing (3.2) we see that

$$\frac{d^m}{ds^m} \frac{1}{(s-1)\varsigma(s)} \bigg|_{s=1} = Y_m \left( -\varepsilon_0, 1!\varepsilon_1, ..., (-1)^m (m-1)!\varepsilon_{m-1} \right)$$

$$= (-1)^m Y_m \left( \varepsilon_0, 1! \varepsilon_1, ..., (m-1)! \varepsilon_{m-1} \right)$$

and hence we see that

$$\frac{d^m}{ds^m} \frac{1}{(s-1)\varsigma(s)} \bigg|_{s=1} = (-1)^m d_m \quad \text{where} \quad d_m > 0$$

Similarly, referring to (2.3.1)

$$\log \frac{1}{2(s-1)\varsigma(s)} = -\sum_{n=1}^{\infty} \frac{b_{n-1}}{n} s^n$$

and (7.2) gives us

$$\frac{1}{2(s-1)\varsigma(s)} = \sum_{n=0}^{\infty} Y_n \left( -b_0, -1! b_1, ..., -(n-1)! b_{n-1} \right) \frac{s^n}{n!}$$

We then see that

$$\frac{d^m}{ds^m} \frac{1}{2(s-1)\varsigma(s)} \bigg|_{s=0} = Y_m \left( -b_0, -1! b_1, ..., -(m-1)! b_{m-1} \right)$$

and employing (2.6) this becomes

$$= Y_m \left( -\mu_0, 1! \mu_1, ..., (-1)^m (m-1)! \mu_{m-1} \right)$$

Using (10.5) we obtain

$$\frac{d^m}{ds^m} \frac{1}{2(s-1)\varsigma(s)} \bigg|_{s=0} = (-1)^m Y_m \left( \mu_0, 1! \mu_1, ..., (m-1)! \mu_{m-1} \right)$$

which has the same sign as $(-1)^m$.

We also see that

$$\frac{d}{ds} \frac{1}{f(s)} = -\frac{f'(s)}{f^2(s)} = -\frac{1}{f(s)} \frac{d}{ds} \log f(s)$$

and applying (10.10) we obtain

$$\frac{d^m}{ds^m} \frac{1}{f(s)} = \frac{1}{f(s)} Y_m \left( g(s), g^{(1)}(s), ..., g^{(m-1)}(s) \right)$$

where $g(s) = -\dfrac{d}{ds}\log f(s)$.

## 9. The $S_2(n)$ constants

We recall (1.16.2)

$$S_2(n) = \sum_{m=1}^{n}\binom{n}{m}\frac{1}{(m-1)!}\frac{d^m}{ds^m}\log[(s-1)\varsigma(s)]\bigg|_{s=1}$$

and we have

$$\frac{d^m}{ds^m}\log[(s-1)\varsigma(s)] = \frac{d^{m-1}}{ds^{m-1}}\frac{f'(s)}{f(s)}$$

where $f(s) = (s-1)\varsigma(s)$. The Leibniz rule for differentiation gives us

$$\frac{d^{m-1}}{ds^{m-1}}f'(s)\frac{1}{f(s)} = \sum_{j=0}^{m-1}\binom{m-1}{j}f^{(m-j)}(s)\frac{d^j}{ds^j}\frac{1}{f(s)}$$

so that

$$\frac{d^m}{ds^m}\log[(s-1)\varsigma(s)]\bigg|_{s=1} = \sum_{j=0}^{m-1}\binom{m-1}{j}(-1)^{m-j-1}(m-j)\gamma_{m-j-1}Y_j\left(\eta_0, 1!\eta_1, ..., (j-1)!\eta_{j-1}\right)$$

We have

$$\frac{d^{n+1}}{ds^{n+1}}[(s-1)\varsigma(s)]\bigg|_{s=1} = (-1)^n(n+1)\gamma_n$$

and

$$\frac{d^{m-j}}{ds^{m-j}}[(s-1)\varsigma(s)]\bigg|_{s=1} = (-1)^{m-j-1}(m-j)\gamma_{m-j-1}$$

and we obtain

$$S_2(n) = \sum_{m=1}^{n}\binom{n}{m}\frac{(-1)^m}{(m-1)!}\sum_{j=0}^{m-1}\binom{m-1}{j}(-1)^{j+1}(m-j)\gamma_{m-j-1}Y_j\left(\eta_0, 1!\eta_1, ..., (j-1)!\eta_{j-1}\right)$$

Using (8.1) we get

$$(-1)^{m-j-1}(m-j)\gamma_{m-j-1} = Y_{m-j}\left(-\eta_0, -1!\eta_1, ..., -(m-j-1)!\eta_{m-j-1}\right)$$

and thus

$$S_2(n) = \sum_{m=1}^{n} \binom{n}{m} \frac{1}{(m-1)!} \sum_{j=0}^{m-1} \binom{m-1}{j} Y_{m-j}\left(-\eta_0, -1!\eta_1, ..., -(m-j-1)!\eta_{m-j-1}\right) Y_j\left(\eta_0, 1!\eta_1, ..., (j-1)!\eta_{j-1}\right)$$

To simplify the notation we write

$$Y_j = Y_j\left(\eta_0, 1!\eta_1, ..., (j-1)!\eta_{j-1}\right)$$

$$Y_j^- = Y_j\left(-\eta_0, -1!\eta_1, ..., -(j-1)!\eta_{j-1}\right)$$

so that

$$S_2(n) = \sum_{m=1}^{n} \binom{n}{m} \frac{1}{(m-1)!} \sum_{j=0}^{m-1} \binom{m-1}{j} Y_{m-j}^- Y_j$$

For example, we see that

$$S_2(1) = Y_0^- Y_1 = \eta_0$$

$$S_2(2) = 2Y_1^- Y_0 + Y_2^- Y_0 + Y_1^- Y_1$$

$$= -2\eta_0 + \eta_0^2 - \eta_1 - \eta_0^2$$

$$= -2\eta_0 - \eta_1$$

It should however be noted that this just adds complexity to the existing problem because we already have the simpler expression (1.17)

$$S_2(n) = -\sum_{m=1}^{n} \binom{n}{m} \eta_{m-1}$$

## 10. Some aspects of the (exponential) complete Bell polynomials

For ease of reference, some aspects of the (exponential) complete Bell polynomials are set out below; these may be useful in the earlier sections of this paper.

As noted, for example, by Kölbig [26] we have

$$(10.1) \qquad \frac{d^r}{dx^r} e^{f(x)} = e^{f(x)} Y_r\left(f^{(1)}(x), f^{(2)}(x), ..., f^{(r)}(x)\right)$$

where the (exponential) complete Bell polynomials may be defined by $Y_0 = 1$ and for $r \geq 1$

(10.2) $\qquad Y_r(x_1, ..., x_r) = \sum_{\pi(r)} \frac{r!}{k_1! \, k_2! \, ... \, k_r!} \left(\frac{x_1}{1!}\right)^{k_1} \left(\frac{x_2}{2!}\right)^{k_2} ... \left(\frac{x_r}{r!}\right)^{k_r}$

where the sum is taken over all partitions $\pi(r)$ of $r$, i.e. over all sets of integers $k_j$ such that

$$k_1 + 2k_2 + 3k_3 + \cdots + rk_r = r$$

The complete Bell polynomials have integer coefficients and the first five are set out below (Comtet [14, p.307])

(10.3) $\qquad\qquad\qquad Y_1(x_1) = x_1$

$$Y_2(x_1, x_2) = x_1^2 + x_2$$

$$Y_3(x_1, x_2, x_3) = x_1^3 + 3x_1 x_2 + x_3$$

$$Y_4(x_1, x_2, x_3, x_4) = x_1^4 + 6x_1^2 x_2 + 4x_1 x_3 + 3x_2^2 + x_4$$

$$Y_5(x_1, x_2, x_3, x_4, x_5) = x_1^5 + 10x_1^3 x_2 + 10x_1^2 x_3 + 15x_1 x_2^2 + 5x_1 x_4 + 10x_2 x_3 + x_5$$

The definition (10.2) immediately implies the following relation

(10.4) $\qquad Y_n(ax_1, a^2 x_2, ..., a^n x_n) = a^n Y_n(x_1, ..., x_n)$

and with $a = 1$ we have

(10.5) $\qquad Y_n(-x_1, x_2, ..., (-1)^n x_n) = (-1)^n Y_n(x_1, ..., x_n)$

Hence if all of the $x_j$ are positive numbers, we see that $Y_n(-x_1, x_2, ..., (-1)^n x_n)$ has the same sign as $(-1)^n$.

By letting $x_j \rightarrow -x_j$ we also note that

(10.6) $\qquad Y_n(x_1, -x_2, ..., (-1)^{n+1} x_n) = (-1)^n Y_n(-x_1, ..., -x_n)$

but no discernible sign pattern emerges here. We may note that

$$Y_n\left(b_1, -1!b_2, ..., (-1)^{n+1}(n-1)!b_n\right) = (-1)^n Y_n\left(-b_1, -1!b_2, ..., -(n-1)!b_n\right)$$

$$= (-1)^n e^{-b_0} \frac{d^n}{dx^n} \frac{1}{h(x)}\bigg|_{x=0}$$

where $h(x)$ is defined by (7.1).

The complete Bell polynomials are also given by the exponential generating function (Comtet [14, p.134])

$$(10.7) \qquad \exp\left(\sum_{j=1}^{\infty} x_j \frac{t^j}{j!}\right) = \sum_{n=0}^{\infty} Y_n(x_1, ..., x_n) \frac{t^n}{n!}$$

Since the exponential function is positive for real arguments, we see that the summation on the right-hand side must also be positive.

Using (10.1) we see that

$$\frac{d^n}{dt^n} \exp\left(\sum_{j=1}^{\infty} x_j \frac{t^j}{j!}\right)\bigg|_{t=0} = Y_n(x_1, ..., x_n)$$

and hence we note that (10.7) is simply the corresponding Maclaurin series.

We note that

$$(10.8) \qquad \sum_{n=0}^{\infty} Y_n(ax_1, ..., ax_n) \frac{t^n}{n!} = \exp\left(\sum_{j=1}^{\infty} ax_j \frac{t^j}{j!}\right) = \exp a\left(\sum_{j=1}^{\infty} x_j \frac{t^j}{j!}\right) = \left[\exp\left(\sum_{j=1}^{\infty} x_j \frac{t^j}{j!}\right)\right]^a$$

and thus we have

$$(10.9) \qquad \left[\sum_{n=0}^{\infty} Y_n(x_1, ..., x_n) \frac{t^n}{n!}\right]^a = \sum_{n=0}^{\infty} Y_n(ax_1, ..., ax_n) \frac{t^n}{n!}$$

Let us now consider a function $f(t)$ which has a Taylor series expansion around $x$: we have

$$e^{f(x+t)} = \exp\left(\sum_{j=0}^{\infty} f^{(j)}(x) \frac{t^j}{j!}\right) = e^{f(x)} \exp\left(\sum_{j=1}^{\infty} f^{(j)}(x) \frac{t^j}{j!}\right)$$

$$= e^{f(x)}\left\{1 + \sum_{n=1}^{\infty} Y_n\left(f^{(1)}(x), f^{(2)}(x), ..., f^{(n)}(x)\right) \frac{t^n}{n!}\right\}$$

We see that

$$\frac{d^m}{dx^m} e^{f(x)} = \frac{\partial^m}{\partial x^m} e^{f(x+t)}\bigg|_{t=0} = \frac{\partial^m}{\partial t^m} e^{f(x+t)}\bigg|_{t=0}$$

and we therefore obtain a derivation of (10.1) above

$$\frac{d^r}{dx^r} e^{f(x)} = e^{f(x)} Y_r\left(f^{(1)}(x), f^{(2)}(x), \ldots, f^{(r)}(x)\right)$$

Suppose that $h'(x) = h(x)g(x)$ and let $f(x) = \log h(x)$. We see that

$$f'(x) = \frac{h'(x)}{h(x)} = g(x)$$

and then using (10.1) above we have

$$(10.10) \qquad \frac{d^r}{dx^r} h(x) = \frac{d^r}{dx^r} e^{\log h(x)} = h(x) Y_r\left(g(x), g^{(1)}(x), \ldots, g^{(r-1)}(x)\right)$$

In particular we have

$$\frac{d^r}{dx^r} h(x)\bigg|_{x=0} = h(0) Y_r\left(g(0), g^{(1)}(0), \ldots, g^{(r-1)}(0)\right)$$

Let us consider the function $h(x)$ with the following Maclaurin expansion

$$(10.11) \qquad \log h(x) = b_0 + \sum_{n=1}^{\infty} \frac{b_n}{n} x^n$$

and we wish to determine the coefficients $a_n$ such that

$$(10.12) \qquad h(x) = \sum_{r=0}^{\infty} a_r x^r$$

By differentiating (10.11) we obtain

$$h'(x) = h(x) \sum_{n=1}^{\infty} b_n x^{n-1} \equiv h(x) g(x)$$

From (10.11) we have

$$\frac{d^r}{dx^r}h(x) = \frac{d^r}{dx^r}e^{\log h(x)} = h(x)Y_r\left(g(x), g^{(1)}(x), ..., g^{(r-1)}(x)\right)$$

and in particular we have

$$\left.\frac{d^r}{dx^r}h(x)\right|_{x=0} = h(0)Y_r\left(g(0), g^{(1)}(0), ..., g^{(r-1)}(0)\right)$$

Using (10.12) the Maclaurin series gives us

$$a_r = \frac{1}{r!}\left.\frac{d^r}{dx^r}h(x)\right|_{x=0}$$

$$a_r = \frac{1}{r!}h(0)Y_r\left(g(0), g^{(1)}(0), ..., g^{(r-1)}(0)\right)$$

We have

$$g^{(j)}(x) = \sum_{n=1}^{\infty}b_n(n-1)(n-2)\cdots(n-j)\,x^{n-1-j}$$

and thus

$$g^{(j)}(0) = j!\,b_{j+1}$$

Therefore we obtain

$$a_r = \frac{1}{r!}e^{b_0}Y_r\left(g(0), g^{(1)}(0), ..., g^{(r-1)}(0)\right)$$

$$= \frac{1}{r!}e^{b_0}Y_r\left(b_1, 1!b_2, ..., (r-1)!b_r\right)$$

Since $\log h(0) = \log a_0 = b_0$ we have

(10.13) $$h(x) = e^{b_0}\sum_{n=0}^{\infty}Y_n\left(b_1, 1!b_2, ..., (n-1)!b_n\right)\frac{x^n}{n!}$$

Then referring to (10.7)

$$\exp\left(\sum_{j=1}^{\infty}x_j\frac{t^j}{j!}\right) = \sum_{n=0}^{\infty}Y_n(x_1, ..., x_n)\frac{t^n}{n!}$$

we see that

$$h(x) = e^{b_0} \exp\left( \sum_{j=1}^{\infty} \frac{b_j}{j} t^j \right)$$

and this is where we started from in (10.12).

$$\log h(x) = b_0 + \sum_{n=1}^{\infty} \frac{b_n}{n} x^n$$

Multiplying (10.12) by $\alpha$ it is easily seen that

$$(10.14) \qquad h^{\alpha}(x) = e^{\alpha b_0} \sum_{n=0}^{\infty} Y_n\left( \alpha b_1, 1! \alpha b_2, ..., (n-1)! \alpha b_n \right) \frac{x^n}{n!}$$

and, in particular, with $\alpha = -1$ we obtain

$$(10.15) \qquad \frac{1}{h(x)} = e^{-b_0} \sum_{n=0}^{\infty} Y_n\left( -b_1, -1! b_2, ..., -(n-1)! b_n \right) \frac{x^n}{n!}$$

Differentiating (10.14) with respect to $\alpha$ would give us an expression for $h^{\alpha}(x) \log h(x)$.

We have the recurrence relation [8, p.415]

$$(10.16) \qquad Y_{r+1}(x_1, ..., x_{r+1}) = \sum_{k=0}^{r} \binom{r}{k} Y_{r-k}(x_1, ..., x_{r-k}) x_{k+1}$$

Using this we find that

$$(r+1)! e^{-b_0} a_{r+1} = \sum_{k=0}^{r} \binom{r}{k} (r-k)! e^{-b_0} a_{r-k} k! b_{k+1}$$

giving us the recurrence relation

$$(r+1) a_{r+1} = \sum_{k=0}^{r} a_{r-k} b_{k+1}$$

or equivalently

$$(10.17) \qquad r a_r = \sum_{m=1}^{r} a_{r-m} b_m$$

Suppose that $h'(x) = h(x) g(x)$ and let $f(x) = \log h(x)$. We see that

$$f'(x) = \frac{h'(x)}{h(x)} = g(x)$$

and then using (10.1) above we have

$$\frac{d^r}{dx^r} h(x) = \frac{d^r}{dx^r} e^{\log h(x)} = h(x) Y_r \left( g(x), g^{(1)}(x), ..., g^{(r-1)}(x) \right)$$

Hence we have the Maclaurin expansion

$$h(x) = h(0) \sum_{r=0}^{\infty} Y_r \left( g(0), g^{(1)}(0), ..., g^{(r-1)}(0) \right) \frac{x^r}{r!}$$

Reference to (10.7) gives us

$$h(x) = h(0) \exp \left( \sum_{j=1}^{\infty} g^{(j-1)}(0) \frac{t^j}{j!} \right)$$

and hence we have

$$\log h(x) = \log h(0) + \sum_{j=1}^{\infty} g^{(j-1)}(0) \frac{t^j}{j!}$$

Therefore, as expected, we obtain

$$f(x) = f(0) + \sum_{j=1}^{\infty} f^{(j)}(0) \frac{t^j}{j!}$$

$\square$

We refer to (7.3)

$$h^{\alpha}(x) = e^{\alpha b_0} \sum_{n=0}^{\infty} Y_n \left( \alpha b_1, 1! \alpha b_2, ..., (n-1)! \alpha b_n \right) \frac{x^n}{n!}$$

which may be written as

$$h^{\alpha}(x) = e^{\alpha b_0} + e^{\alpha b_0} \sum_{n=1}^{\infty} Y_n \left( \alpha b_1, 1! \alpha b_2, ..., (n-1)! \alpha b_n \right) \frac{x^n}{n!}$$

We have [8, p.412]

$$Y_n(x_1,...,x_n) = \sum_{k=1}^{n} B_{n,k}(x_1,...,x_{n-k+1})$$

so that

$$Y_n(\alpha x_1,...,\alpha x_n) = \sum_{k=1}^{n} B_{n,k}(\alpha x_1,...,\alpha x_{n-k+1})$$

We also have [8, p.412]

$$B_{n,k}(\alpha x_1,...,\alpha x_{n-k+1}) = \alpha^k B_{n,k}(x_1,...,x_{n-k+1})$$

and hence we see that

$$(10.18) \quad Y_n(\alpha x_1,...,\alpha x_n) = \sum_{k=1}^{n} \alpha^k B_{n,k}(x_1,...,x_{n-k+1})$$

Therefore we obtain

$$(10.19) \quad h^\alpha(x) = e^{\alpha b_0} + e^{\alpha b_0} \sum_{n=1}^{\infty} \sum_{k=1}^{n} \alpha^k B_{n,k}\left(b_1, 1!b_2,...,(n-k)!b_{n-k+1}\right)\frac{x^n}{n!}$$

Differentiating (10.19) with respect to $\alpha$ gives us

$$h^\alpha(x)\log h(x) = b_0 e^{\alpha b_0} + b_0 e^{\alpha b_0} \sum_{n=1}^{\infty} \sum_{k=1}^{n} \alpha^k B_{n,k}\left(b_1, 1!b_2,...,(n-k)!b_{n-k+1}\right)\frac{x^n}{n!}$$

$$+ e^{\alpha b_0} \sum_{n=1}^{\infty} \sum_{k=1}^{n} k\alpha^{k-1} B_{n,k}\left(b_1, 1!b_2,...,(n-k)!b_{n-k+1}\right)\frac{x^n}{n!}$$

and letting $\alpha = 0$ results in

$$\log h(x) = b_0 + \sum_{n=1}^{\infty} B_{n,1}\left(b_1, 1!b_2,...,(n-1)!b_n\right)\frac{x^n}{n!}$$

Comparing this with (7.1) shows that

$$B_{n,1}\left(b_1, 1!b_2,...,(n-1)!b_n\right) = (n-1)!b_n$$

A further differentiation gives us

$$h^{\alpha}(x)\log^2 h(x) = b_0^2 e^{\alpha b_0} + b_0^2 e^{\alpha b_0} \sum_{n=1}^{\infty}\sum_{k=1}^{n} \alpha^k B_{n,k}\left(b_1, 1!b_2, ..., (n-k)!b_{n-k+1}\right)\frac{x^n}{n!}$$

$$+ 2b_0 e^{\alpha b_0} \sum_{n=1}^{\infty}\sum_{k=1}^{n} k\alpha^{k-1} B_{n,k}\left(b_1, 1!b_2, ..., (n-k)!b_{n-k+1}\right)\frac{x^n}{n!}$$

$$+ e^{\alpha b_0} \sum_{n=1}^{\infty}\sum_{k=1}^{n} k(k-1)\alpha^{k-2} B_{n,k}\left(b_1, 1!b_2, ..., (n-k)!b_{n-k+1}\right)\frac{x^n}{n!}$$

and letting $\alpha = 0$ results in

$$\log^2 h(x) = b_0^2 + 2b_0 \sum_{n=1}^{\infty} B_{n,1}\left(b_1, 1!b_2, ..., (n-1)!b_n\right)\frac{x^n}{n!} + 2\sum_{n=1}^{\infty} B_{n,2}\left(b_1, 1!b_2, ..., (n-2)!b_{n-1}\right)\frac{x^n}{n!}$$

or, since $B_{1,2} = 0$, we have

$$(10.20) \qquad \log^2 h(x) = b_0^2 + 2b_0\left[\log h(x) - b_0\right] + 2\sum_{n=2}^{\infty} B_{n,2}\left(b_1, 1!b_2, ..., (n-2)!b_{n-1}\right)\frac{x^n}{n!}$$

Differentiating (10.18) gives us

$$\frac{\partial^m}{\partial\alpha^m} Y_n(\alpha x_1, ..., \alpha x_n) = \sum_{k=1}^{n} k(k-1)...(k-m+1)\alpha^{k-m} B_{n,k}(x_1, ..., x_{n-k+1})$$

and we obtain

$$\left.\frac{\partial^m}{\partial\alpha^m} Y_n(\alpha x_1, ..., \alpha x_n)\right|_{\alpha=0} = m! B_{n,m}(x_1, ..., x_{n-m+1})$$

$\square$

Further applications of the (exponential) complete Bell polynomials are contained in [18]. For example, it is shown in [18] how may they be employed to prove that $\Gamma^{(n)}(x)$ has the same sign as $(-1)^n$ for all $x$ in the interval $(0, \alpha)$ where $\alpha > 0$ is the unique solution of $\psi(\alpha) = 0$.

Donal F. Connon
Elmhurst
Dundle Road
Matfield
Kent TN12 7HD
dconnon@btopenworld.com